\documentclass[11pt]{article} 

\usepackage{latexsym} 
\usepackage{amsfonts,amssymb} 
\usepackage{amsmath} 
\usepackage{amsthm} 
\usepackage{mathrsfs} 
\usepackage{dsfont} 
\usepackage{bbold} 
\usepackage[english]{babel} 
\usepackage{caption} 
\usepackage{epsfig} 
\usepackage{float} 
\usepackage{color} 
\usepackage{subfigure} 
\usepackage{psfrag} 
\usepackage{graphicx} 
\usepackage{epsfig} 
\usepackage{amsbsy} 
\usepackage{hyperref}

\newtheorem{theorem}{Theorem} 
\newtheorem*{prop*}{Theorem} 
\newtheorem{thm}[theorem]{Theorem} 
\newtheorem{coro}[theorem]{Corollary} 
\newtheorem{defi}[theorem]{Definition} 
\newtheorem{lemma}[theorem]{Lemma} 
\newtheorem{prop}[theorem]{Proposition} 
\newtheorem{rmk}[theorem]{Remark}

\newcommand{\zerarcounters}{\setcounter{equation}{0}\setcounter{theorem}{0}} 

\newcommand{\beq}{\begin{equation}}
\newcommand{\eeq}{\end{equation}}

\newcommand{\ZZZ}{\mathds{Z}} 
\newcommand{\CCC}{\mathds{C}} 
\newcommand{\NNN}{\mathds{N}} 
 
\newcommand{\RRR}{\mathds{R}} 
\newcommand{\TTT}{\mathds{T}} 
\newcommand{\uno}{\mathds{1}} 

\newcommand{\calA}{{\mathcal A}} 
\newcommand{\BB}{{\mathcal B}} 
 
\newcommand{\DD}{{\mathcal D}} 
\newcommand{\calE}{{\mathcal E}} 
 
\newcommand{\calG}{{\mathcal G}} 
\newcommand{\calH}{{\mathcal H}}

\newcommand{\LL}{{\mathcal L}} 
\newcommand{\MM}{{\mathcal M}} 
 
\newcommand{\calO}{{\mathcal O}} 
\newcommand{\calP}{{\mathcal P}} 
 
\newcommand{\RR}{{\mathcal R}} 
 
\newcommand{\TT}{{\mathcal T}} 
\newcommand{\calU}{{\mathcal U}}

\newcommand{\gotf}{{\mathfrak f}}

\newcommand{\ol}{\overline} 

\newcommand{\Fullbox}{{\rule{2.0mm}{2.0mm}}} 
 
\newcommand{\EP}{\hfill\Fullbox\vspace{0.2cm}} 
\newcommand{\prova}{\noindent{\it Proof. }} 
\newcommand{\io}{\infty} 
\newcommand{\e}{\varepsilon} 
 
\newcommand{\al}{\alpha} 
\newcommand{\de}{\delta} 
\newcommand{\be}{\beta} 
 
\newcommand{\m}{\mu} 
\newcommand{\x}{\xi}

\newcommand{\g}{\gamma} 
\newcommand{\om}{\omega} 
\newcommand{\h}{\eta} 
\newcommand{\ze}{\zeta} 
\newcommand{\la}{\lambda} 
 
\newcommand{\f}{\varphi} 
\newcommand{\s}{\sigma} 
 
\newcommand{\del}{\partial}

\newcommand{\avg}[1]{\langle #1 \rangle}

\newcommand{\ii}{{\rm i}}

\newcommand{\tc}{{\mathtt{c}}}

\newcommand{\tC}{{\mathtt{C}}}

\newcommand{\cF}{{\mathcal F}}

\newcommand{\cH}{{\mathcal H}}

\newcommand{\cL}{{\mathcal L}}

\def\tilde#1{\widetilde{#1}}

\def\ins#1#2#3{\vbox to0pt{\kern-#2 \hbox{\kern#1 #3}\vss}\nointerlineskip} 

\newcommand{\zia}{\eta}
\newcommand{\N}{{\mathds N}}

\newcommand{\C}{{\mathds C}}
\newcommand{\R}{{\mathds R}}

\newcommand{\T}{{\mathds T}}
\newcommand{\Z}{{\mathds Z}}

\newcommand{\Oo}{\Omega}
\newcommand{\Dc}{{\mathtt{ D}_\g}}
\newcommand{\Dcb}{{\mathtt{ D}_{\ol{\g}}}}
\newcommand{\im}{{\rm i}}
\newcommand{\jap}[1]{\langle #1 \rangle}
\newcommand{\set}[1]{{\left\{#1\right\}}}
\newcommand{\pa}[1]{{\left(#1\right)}}

\newcommand{\Id}{\mathtt{ Id}}
\begin{document}
 
%%%%%%%%%%%%%%%%%%%%%%%%%%%%%%%%%%%%%%%%%%%%%%%%%
\title{\bf Almost-periodic Response Solutions for a forced quasi-linear Airy equation}
%%%%%%%%%%%%%%%%%%%%%%%%%%%%%%%%%%%%%%%%%%%%%%%%%
 
 \author{\bf Livia Corsi$^\dag$, Riccardo Montalto$^*$, Michela Procesi$^\dag$
\\
\small ${}^*$ Universit\`a degli Studi di Milano; riccardo.montalto@unimi.it \\
\small
${}^\dag$ Universit\`a di Roma Tre; lcorsi@mat.uniroma3.it, procesi@mat.uniroma3.it}
\date{} 
 
\maketitle

\begin{abstract} 
We prove the existence of almost-periodic solutions for quasi-linear perturbations
of the Airy equation. This is the first result about the existence of this type of solutions
for a quasi-linear PDE. The solutions turn out to be analytic in time and space.
To prove our result we use a Craig-Wayne approach combined with a KAM reducibility scheme
and pseudo-differential calculus on $\TTT^\io$.

\smallskip


\smallskip

\end{abstract}   
  
\tableofcontents  

%%%%%%%%%%%%%%%%%%%%%%%%%%%%%%%%%%%%%%%%%%%%%%%%%
%%%%%%%%%%%%%%%%%%%%%%%%%%%%%%%%%%%%%%%%%%%%%%%%%
\zerarcounters
\section{Introduction}
%%%%%%%%%%%%%%%%%%%%%%%%%%%%%%%%%%%%%%%%%%%%%%%%%
%%%%%%%%%%%%%%%%%%%%%%%%%%%%%%%%%%%%%%%%%%%%%%%%%

In this paper we study response solutions for almost-periodically forced quasilinear PDEs
close to an elliptic fixed point.

The problem of response solutions for PDEs has been widely studied in many contexts,
starting from the papers \cite{Rab1,Rab2}, where the Author considers a periodically forced PDE
with dissipation. In the presence of dissipation, of course there is no small divisors problem.
However as soon as the dissipation is removed, small divisors appear even in the easiest
possible case of a periodic forcing when the spacial variable is one dimensional.

The first results of this type in absence of dissipation were obtained by means of a KAM
approach \cite{K,K2,K3,W,P,KP}. However, a more functional approach, via a combination of a 
Ljapunov-Schmidt reduction and a Newton scheme, in the spirit of \cite{Rab1,Rab2},
was proposed by Craig-Wayne \cite{CW}, and then generalized in many ways by Bourgain;
see for instance \cite{B1,B2,B3} to mention a few.
All the results mentioned above concern semi-linear PDEs and the forcing is quasi-periodic.

In more recent times, the Craig-Wayne-Bourgain approach has been fruitfully used and
generalized in order to cover quasi-linear and fully nonlinear PDEs, again in the quasi-periodic case;
see for instance \cite{BBM,FP,CM,BBHM}.

Regarding the almost-periodic case, most of the classical results are obtained via a KAM-like approach;
see for instance \cite{CP,Po,Bjfa}. A notable exception is \cite{Bgafa}, where the Craig-Wayne-Bourgain
method is used.
More recently there have been results such as \cite{RLZ,RL1,Liu}, which use a KAM approach.
We mention also \cite{XG,Y,BMP1,BMP} which however are tailored for an autonomous PDE.

All the aforementioned results, concern semi-linear PDEs, with no derivative in the nonlinearity.
Moreover they require a very strong analyticity condition on the forcing term.
Indeed the difficulty of proving the existence of almost-periodic response solution is strongly
related to the regularity of the forcing, since one can see an almost periodic function 
as the limit of quasi-periodic ones with an increasing number
of frequencies.
 If such limit is reached sufficiently fast, the most direct strategy would be to
 iteratively find approximate quasi-periodic response solutions and then take the limit.
 This is the overall strategy of \cite{Po} and \cite{RLZ,RL1,Liu}.
 However this procedure works if one considers a sufficiently
regular forcing term and a bounded nonlinearity, but becomes very delicate in the case of unbounded 
 nonlinearities.

In the present paper we study the existence of almost-periodic response solutions, for a 
quasi-linear PDE on $\TTT$. To the best of our knowledge this is the first result of this type.

Specifically we consider a quasi-linear Airy equation 
\begin{equation}\label{kdv}
\partial_t  u +  \partial_{xxx} u   +  Q(u, u_x, u_{xx}, u_{xxx}) + 
 \gotf( t, x) =0,\qquad x\in\TTT := (\RRR /(2 \pi \ZZZ))
\end{equation}
where $Q$ is a Hamiltonian, 
quadratic nonlinearity and $\gotf$ is an analytic forcing term
 with zero average w.r.t. $x$.  We assume $\gotf$ to be ``almost-periodic'' with frequency
 $\omega\in\ell^{\io}$, in the sense of Definition \ref{suca}.
 
We mention that in the context of reducibility of linear PDEs a problem of this
kind has been solved in \cite{MP}. Our aim is to provide a link between the linear techniques
of \cite{MP} and the nonlinear Craig-Wayne-Bourgain method. Note that such a link is
nontrivial, and requires a delicate handling; see below.
 
 The overall setting we use is the one of \cite{BBM}. However their strategy is taylored
 for Sobolev regularity; the quasi-periodic analytic case has been covered in \cite{CFP}. Unfortunately
 the ideas of \cite{CFP} cannot be directly applied in the almost-periodic case.
Roughly, it is well known that the regularity and the small-divisor problem conflict. Thus, in
the almost-periodic case one expect this issue to be even more dramatic.
Specifically, we were not able to define a ``Sobolev'' norm for almost-periodic functions,
satisfying the interpolation estimates needed in the Nash-Moser scheme;
 this is why we cannot  use the theorem of \cite{CFP}.
 
 Let us now present our main result in a more detailed way.
 
First of all we note that \eqref{kdv} is an Hamiltonian PDE whose
 Hamiltonian  is given by 
\begin{equation}\label{Hamiltoniana}
H(u) := \frac{1}{2}\int_\TTT  u_x^2 d x -  \frac{1}{6} \int_\TTT G(u, u_x)\, d x
 - \int_\T F( t, x)  u  d x ,\qquad
\gotf(t,x) = \del_x F(t,x)
\end{equation}
where $G(u, u_x)$ is a cubic Hamiltonian density of the form 
\begin{equation}\label{forma densita hamiltoniana quadratica}
G(u, u_x) := \mathtt c_3 u_x^3 + \mathtt c_2 u u_x^2 + \mathtt c_1 u^2 u_x + 
\mathtt c_0 u^3, \quad \mathtt c_0, \ldots, \mathtt c_3 \in \R\,
\end{equation}
and the symplectic structure is given by $J=\del_x$.
The Hamiltonian nonlinearity $Q(u, \ldots, u_{xxx})$ is therefore given by 
\begin{equation}\label{nonlin quadratica Q fQ}
Q(u, u_x, u_{xx}, u_{xxx})  = \partial_{xx}(\partial_{u_x} G(u, u_x)) - \partial_x (\partial_u G(u, u_x)) 
\end{equation}
and the Hamilton equations are 
$$
\partial_t u = \partial_x \nabla_{ u} H(u)\,. 
$$

We look for an almost-periodic solution to \eqref{kdv} with frequency $\om$ in the sense below.

For $\zia>0$, define the set of infinite integer vectors with {\it finite support}  as
\begin{equation}\label{Z inf *} 
\Z^\infty_* :=   \Big\{ \ell \in \Z^\N : |\ell|_\zia := \sum_{i\in \N}  i ^\zia |\ell_i| < \infty \Big\}. 
\end{equation}
Note that $\ell_i \neq 0$ only for finitely many indices $i \in \N$. In particular $\ZZZ^\io_*$
does not depend on  $\h$.

\begin{defi}\label{suca}
	Given $\omega\in [1,2]^{\NNN}$ with rationally independent components\footnote{We say
	that $\om$ has rationally independent components if for any $N>0$ and any $k\in\ZZZ^N$ one has
	$\sum_{i=1}^N \om_i k_i \ne0$.}
	and a Banach space $(X,|\cdot|_X)$,  we say that $F(t):\R\to X$ 
	is almost-periodic in time with frequency $\omega$ and analytic in the strip $\s>0$ 
	if we may write it in totally   convergent Fourier series
	\[
	\begin{aligned}
	F(t)= &\sum_{ \ell \in \Z^\infty_*}  F(\ell) e^{\im \ell \cdot \omega t }
	\quad \text{such that} \quad  F(\ell) \in X\,,\;\forall \ell \in \Z^\infty_* \quad 
	\\
	&\text{and}
	\quad
	|F|_\s:=
	\sum_{ \ell \in \Z^\infty_*} |  F(\ell) |_X e^{\s |\ell|_\zia } <\infty.
	\end{aligned}
	\]
\end{defi}

We shall be particularly interested in almost-periodic functions where $X=\cH_0(\T_\s)$ 
\[
\cH_0(\T_\s):= \Big\{ u= \sum_{j \in\Z\setminus\{0\}} u_j e^{\im j x}\,,\;  u_j=\bar u_{-j} \in 
\C\,:\quad |u|_{\cH(\T_\s)}:= \sum_{j \in\Z \setminus \{ 0 \}}| u_j| e^{\s |j| } <\infty \Big\}
\]
is the space of analytic, real on real functions $\T_s\to \C$ with zero-average, where   
$\T_s := \{ \f \in \C : {\rm Re}(\f) \in \T, \ |{\rm Im}(\f)| \leq s \}$ is the thickened torus. We recall that a function $u : \T_s \to \C$ is real on real if for any $x \in \T$, $u(x) \in \R$. 
\smallskip

Of course we need some kind of Diophantine condition on $\om$. We give the following,
taken from \cite{Bjfa,MP}.

\begin{defi}\label{diofantino} Given $\gamma \in (0, 1)$, we  denote by $\Dc$ the set of 
{\it Diophantine} frequencies  
	\begin{equation}\label{diofantinoBIS}
	\Dc:=\set{\omega\in [1,2]^\N \,:\;	|\omega\cdot \ell|> \gamma 
	\prod_{i\in \N}\frac{1}{(1+|\ell_i|^{2} {i}^{2})}\,,
		\quad \forall \ell\in \Z^\io_*\setminus\{0\}}.
	\end{equation} 
	\end{defi}

We are now ready to state our main result.

\begin{thm}[Main Theorem]\label{erteo}
Fix $\ol{\g}$. 
	Assume  that $\gotf$ in \eqref{kdv} is almost-periodic in time and
	 analytic in a strip $S$ (both in time and space). Fix $\ol{s}<S$.
	  If $\gotf$ has an appropriately small norm depending on $S-\ol{s}$, namely
	\begin{equation}\label{certo}
	|\gotf|_S := \sum_{ \ell \in \Z^\infty_*} |  \gotf(\ell) |_{\cH_0(\T_S)} e^{S |\ell|_\zia }  \le \epsilon(S-\ol{s})\ll 1,
	\qquad \epsilon(0)=0,
	\end{equation}
	then there is a Cantor-like set $\calO^{(\io)}\subseteq \Dcb$
	with positive Lebesgue measure, and for all $\om\in\calO^{(\io)}$
	a solution to \eqref{kdv}  which is almost-periodic in time with frequency $\om$ and
	 analytic in a strip $\ol{s}$ (both in time and space).
	\end{thm}
\begin{rmk}
	Of course the same result holds verbatim if we replace the quadratic polynomial $Q$ by a polinomial of arbitrary degree.
	We could also assume that
	the coefficients $\tc_j$ appearing in \eqref{nonlin quadratica Q fQ} depend on $x$ and $\omega t$. In that case Theorem \ref{erteo}  holds provided we further require a condition of the type
$
\sup_{j} |\partial_x^2 \tc_j|_S \le C$. Actually one could also take $Q$ to be an analytic
function with a zero of order two. However this leads to a number of long and non particularly
enlightening calculations.
\end{rmk}

To prove Theorem \ref{erteo} we proceed as follows.
First of all we regard \eqref{kdv} as a functional Implicit Function Problem on some appropriate space
of functions defined on an infinite dimensional torus; see Definition \ref{anali} below.
Then in Section \ref{primo} we prove an iterative ``Nash-Moser-KAM'' scheme to produce the solution
of such Implicit Function Problem. 
It is well known that an iterative rapidly converging scheme heavily relies on a careful control
on the invertibility of the linearized operator at any approximate solution. Of course, in the case of a quasi-linear PDE this
amounts to study an unbounded non-constant coefficients operator. To deal with this
problem, at each step
we introduce a change of variables $T_n$ which diagonalizes the highest order terms
of the linearized operator. An interesting feature is that $T_n$ preserves the PDE structure.
As in \cite{CFP} and differently from the classical papers, at each
step we apply the change of variables $T_n$ to the whole nonlinear operator.
This is not a merely technical issue.
Indeed, the norms we use are strongly coordinate-depending, and the change of variable $T_n$ that
we need to apply are not close-to-identity, in the sense that $T_n - \Id$ is not a bounded operator small in size.

In Section \ref{dimo} we show how to construct the change of variables $T_n$ satisfying the properties
above.
Then in order to prove the invertibility of the linearized operator after the change of
variables $T_n$ is applied, one needs to perform a reducibility scheme: this is done in
Section \ref{eppoi}. For a more detailed description of the technical aspects see Remark \ref{dellintro}.

\noindent
{\sc Acknowledgements.} Riccardo Montalto is supported by INDAM-GNFM. 

%%%%%%%%%%%%%%%%%%%%%%%%%%%%%%%%%%%%%%%%%%%%%%%%%
%%%%%%%%%%%%%%%%%%%%%%%%%%%%%%%%%%%%%%%%%%%%%%%%%

\zerarcounters
\section{Functional setting}\label{noianormale}

As it is habitual in the theory of quasi-periodic functions we shall study  almost periodic functions 
in the context of analytic functions on an infinite dimensional torus. To this purpose, for $\zia,s>0$, 
we define the  {\it thickened} infinite dimensional torus $\T^\infty_s$ as 
\[
\f=(\f_i)_{i\in \N}\,,\quad \f_i\in\C\,:\; {\rm Re}(\f_i)\in \T\;,\;|{\rm Im}(\f_i)|\le s \jap{i}^\zia\,. 
\]
Given a Banach space $(X, | \cdot |_X )$ we consider the space $\cF$ of 
pointwise absolutely convergent formal Fourier series $\T^\infty_s \to X$
\begin{equation}
\label{fourier}
u(\f) = \sum_{\ell \in \Z^\infty_*}  u(\ell) e^{\im \ell \cdot \f }\,,\quad  u(\ell) \in X
\end{equation}
and define the analytic functions as follows.
\begin{defi}\label{anali}
	Given a Banach space $(X, | \cdot |_X )$ and $s > 0$, we define the space of analytic 
	functions $\T^\infty_s \to X$ as the  subspace 
	$$
	{\mathcal H}( \T^\infty_s, X) := \Big\{ u(\f) = \sum_{\ell \in \Z^\infty_*}  u(\ell) 
	e^{\im \ell \cdot \f}\in \cF \; \;:\quad | u |_{s} := \sum_{\ell \in \Z^\infty_*} 
	e^{s |\ell|_\zia} |  u(\ell)|_X < \infty  \Big\}\,. 
	$$
\end{defi}
We denote by $\calH_s$ the subspace of ${\mathcal H}(\T^\infty_s, {\mathcal H}_0(\T_s))$ of the functions which are real on real. Moreover, we denote by $ {\mathcal H}(\T^\infty_{s}\times\TTT_s)$, the space of analytic functions $\T^\infty_s \times \T_s \to \C$ which are real on real. The space ${\cal H}_s$ can be identified with the subspace of zero-average functions of $ {\mathcal H}(\T^\infty_{s}\times\TTT_s)$. Indeed if $u \in  {\cal H}_s$, then   
\[
\begin{aligned}
& u= \sum_{ \ell \in \Z^\infty_*} u(\ell,x)e^{\im  \ell \cdot \f}= 
\sum_{ (\ell,j )\in \Z^\infty_*\times \Z\setminus\{0\}} u_j(\ell)e^{\im  \ell \cdot \f +\im j x}, \\
& \text{with} \quad u_j(\ell) = \overline{u_{- j}(- \ell)}
\end{aligned}
\]
For any  $u\in  {\mathcal H}(\T^\infty_{s}\times\TTT_s)$ let us denote
\begin{equation}\label{ricca' basta co sti label lunghissimi porchiddi}
(\pi_0 u)(\f,x) := \avg{u(\f,\cdot)}_x := \frac{1}{2 \pi}\int_\T u(\f, x)\, d x,\qquad \pi_0^\perp := \uno-\pi_0\,.
\end{equation}

Throughout the algorithm we shall need to control the Lipschitz variation w.r.t. $\om$ of functions 
in some $\cH(\T^\io_s,X)$, which are defined for $\om$ in some  Cantor set. Thus, for 
$\calO\subset\calO^{(0)}$
we introduce the following norm.

\medskip

\noindent
{\bf Parameter dependence.} Let $Y$ be a Banach space and $\gamma \in (0, 1)$. If 
$f : \Oo \to Y$, $\Oo \subseteq  [1, 2]^\N$ is a Lipschitz function we define 
\begin{equation}\label{normalip}
\begin{aligned}
& | f |_Y^{\rm sup} := \sup_{\omega \in \Oo } | f(\omega)|_Y, \quad | f |_Y^{\rm lip} := \sup_{\begin{subarray}{c}
	\omega_1, \omega_2 \in \Oo \\
	\omega_1 \neq \omega_2
	\end{subarray}} \frac{| f(\omega_1) - f(\omega_2) |_Y}{| \omega_1 - \omega_2|_\infty}\,, \\
& | f |_Y^\Oo := | f |_Y^{\rm sup} + \gamma | f |_Y^{\rm lip}\,. 
\end{aligned}
\end{equation}
If $Y = \calH_s$ we simply write $| \cdot |_\sigma^{\rm sup}$, $| \cdot |_\sigma^{\rm lip}$, 
$| \cdot |_\sigma^{\Oo}$. If $Y$ is a finite dimensional space, we write $| \cdot |^{\rm sup}$, 
$| \cdot |^{\rm lip}$, $| \cdot |^{\Oo}$.

%\begin{equation}\label{normalip}
%|g|_s^{\calO} := \sup_{\om\in\calO}|g(\cdot,\om)|_s + \g \sup_{\om_1\ne\om_2\in\calO}
%\frac{|g(\cdot,\om_1)-g(\cdot,\om_2)|_s}{|\om_1-\om_2|_{\io}} 
%\end{equation}
%and denote
%\begin{equation}\label{analilip}
%\calH_{s}^{\calO}:= \{ g=g(\cdot,\om)= \sum_{\ell\in\ZZZ^\io_*,j\in\ZZZ} 
%g_{j}(\ell)e^{\im (\f\cdot\ell +jx)}\;:\; |g|_s^{\calO}<\io\}.
%\end{equation}

\medskip

\noindent
{\bf Linear operators.} For any $\sigma > 0$, $m \in \R$ we define the class of linear 
operators of order $m$ (densely defined on $L^2(\T)$) ${\mathcal B}^{\sigma, m}$ as 
\begin{equation}\label{definizione classe cal B sigma}
\begin{aligned}
& {\mathcal B}^{\sigma, m} := \Big\{ {\mathcal R} : L^2(\T) \to L^2(\T) : 
\| {\mathcal R}\|_{{\mathcal B}^{\sigma, m}} < \infty \Big\} \quad \text{where} \\
& \quad \| {\mathcal R}\|_{{\mathcal B}^{\sigma, m}}
 := \sup_{j' \in \ZZZ\setminus\{0\}} \sum_{j \in \ZZZ\setminus\{0\}} 
e^{\sigma|j - j'| } |R_j^{j'}| \langle j' \rangle^{- m} \,. 
\end{aligned}
\end{equation}
and for $\TT \in {\mathcal H}(\T^\infty_\sigma, {\mathcal B}^{\sigma, m})$ we set
\begin{equation}\label{definizione norma del decay}
\|{\mathcal T}\|_{\sigma, m} := \sum_{\ell \in \Z^\infty_*} e^{\sigma |\ell|_\zia}
 \| {\mathcal T}(\ell)\|_{{\mathcal B}^{\sigma, m}}\,. 
\end{equation}

In particular we shall denote by $\| \cdot \|_{{\sigma, m}}^{\Omega}$ the corresponding
Lipshitz norm. Moreover if $m=0$ we shall drop it, and write simply 
$\| \cdot \|_{{\sigma}}$ or $\| \cdot \|_{{\sigma}}^{\Omega}$.

%%%%%%%%%%%%%%%%%%%%%%%%%%%%%%%%%%%%%%%%%%%%%%%%%
%%%%%%%%%%%%%%%%%%%%%%%%%%%%%%%%%%%%%%%%%%%%%%%%%

\zerarcounters
\section{The iterative scheme}\label{primo}

Let us rewrite \eqref{kdv} as
\begin{equation}\label{ezero}
F_0(u) =0
\end{equation}
where
\begin{equation}\label{Fzero}
F_0(u):= (\om\cdot\del_\f + \del_{xxx})u + Q(u,u_x,u_{xx},u_{xxx}) + f(\f,x)
\end{equation}
where we $\gotf(t,x)=f(\om t,x)$ and , as custumary the unknown $u$ is a function of 
$(\f,x)\in \TTT^\io \times \TTT$.

We introduce the (Taylor) notation
\begin{equation}\label{tayzero}
\begin{aligned}
& L_0 := (\om\cdot\del_\f + \del_{xxx})=F'_0(0),\qquad f_0=F_0(0) =  f(\f,x),\qquad \\
& Q_0(u) =Q(u,u_x,u_{xx},u_{xxx})\stackrel{\eqref{nonlin quadratica Q fQ}}{=}  
\partial_{xx}\Big(3 \mathtt c_3 u_x^2 + 2 \mathtt c_2 u u_x + \mathtt c_1 u^2  \Big) \\
& \qquad\;\; - \partial_x (\mathtt c_2 u_x^2 + 2 \mathtt c_1 u u_x + 3 \mathtt c_0 u^2 )
\end{aligned}
\end{equation}
so that \eqref{zero} reads
\[
f_0 + L_0 u + Q_0(u) =0.
\]

%Note that $F_0$ is a quadratic operator, mapping $\calH_s$ onto $\calH_{\tilde s}$ where
%$\tilde s< \min\{S,s\}$. 

Note that $Q_0$ is of the form
\begin{equation}\label{q0}
Q_0(u)= \sum_{0\le i\le 2,\,0\le j\le 3 \atop 0\le i+j\le 4} q^{(0)}_{i,j}  (\del_x^i u)(\del_x^j u)
\end{equation}
with the coefficients $q^{(0)}_{i,j}  $ satisfying
\begin{equation}\label{identica0}
\sum_{0\le i\le 2,\,0\le j\le 3 \atop 0\le i+j\le 4} |q^{(0)}_{i,j}|\leq C\,, 
\end{equation}
where the constant $C$ depends clearly on $|\mathtt c_0|, \ldots, |\mathtt c_3|$. 
In particular, this implies that for all $u\in\calH_s$ one has the following.

\begin{itemize}
	
	\item[Q1.] $|Q_0(u)|_{s-\s} \lesssim  \s^{-4}|u|_s^2$
	
	\item[Q2.] $|Q_0'(u)[h]|_{s-\s}\lesssim \s^{-4}|u|_s |h|_s$
	%	\item[Q3.] $|Q_0''(u)[h,k]|_{s-\s}\lesssim \s^{-4}|h|_s|k|_s$
	
\end{itemize}

We now fix the constants
\begin{equation}\label{costanti}
\begin{aligned}
& \mu > \max\{ 1 , \frac{1}{\eta}\} \,, \\
&\g_0<\frac{1}{2}\ol{\g},\qquad \g_n:=(1-2^{-n})\g_{n-1}\,,\quad n\ge1\\
&\s_{-1}:={\frac{1}{8}\min\{(S-\ol{s}),1\}}\,,\qquad
\s_{n-1}=\frac{6\s_{-1}}{\pi^2n^2},\quad n\ge1\,,\\
&s_0=S-\s_{-1}\,,\qquad s_n=s_{n-1}-6\s_{n-1},\quad n\ge1, \\
&\e_n:= \e_0 e^{-\chi^n},\quad \chi=\frac{3}{2}\,,
\end{aligned}
\end{equation}
where $\e_0$ is such that
\begin{equation}\label{quantifico}
e^{\tC_0\s_{-1}^{-\mu }}|f|_S=e^{\tC_0\s_{-1}^{-\mu}}|f_0|_S\ll \e_0.
\end{equation}

Introduce
\begin{equation}\label{def eta (ell)}
\mathtt d(\ell) := \prod_{i\in \N}(1+|\ell_i|^{5} \jap{i}^{5}), \quad \forall \ell \in \Z^\infty_*\,.
\end{equation}

We also set $\calO^{(-1)}:= \mathtt D_{\ol{\g}}$ and
\begin{equation}\label{calo0}
	\calO^{(0)}:=\set{\om\in\Dcb:\quad |\om\cdot\ell + j^3|\ge 
		\frac{\gamma_0}{{\mathtt d}(\ell)}\,,
		\quad \forall \ell\in \Z^\io_*\,,\quad  j\in \N \,,\; (\ell,j)\ne (0,0) }.
	\end{equation}

%Note that $\calO^{(0)}$ has positive measure in $[-1,1]^{\NNN}$ endowed with the
%probability product measure; see Section \ref{misure}.

\begin{prop}\label{ittero}
There exists $\tau,\tau_1,\tau_2,\tau_3, \tC,\epsilon_0$ (pure numbers) such that for
\begin{equation}
\label{piccolo0}
\e_0 \le \s_0^{\tau} e^{-\tC\s_0^{-\mu}} \epsilon_0\,,
\end{equation}
 for all $n\ge1$ the following hold.

\begin{enumerate}

\item
There exist a sequence of Cantor sets $\calO^{(n)}\subseteq \calO^{(n-1)}$, $n\ge1$ such that
\begin{equation}\label{alvolo}
{\mathbb P} (\calO^{(n-1)}\setminus \calO^{(n)}) \lesssim  \frac{\g_0}{n^2}\,.
\end{equation}

\item
For $n\ge 1$, there exists a sequence of linear, invertible, bounded and symplectic changes of variables 
defined for $\om \in \calO^{(n-1)}$, of the form
\begin{equation}\label{tika}
T_{n} v(\f,x)=  (1+\xi^{(n)}_x)v( \f +\omega \be^{(n)}(\f), x+ \xi^{(n)}(\f,x) + p^{(n)}(\f))\,
\end{equation}
satisfying
	\begin{equation}\label{male}
	| \xi^{(n)}|_{s_{n-1}- \s_{n-1}}^{\calO^{(n-1)}},|\be^{(n)}|_{s_{n-1}- \s_{n-1}}^{\calO^{(n-1)}},
	|p^{(n)}|_{s_{n-1}- \s_{n-1}}^{\calO^{(n-1)}} \lesssim \s_{n-1}^{-\tau_1} \e_{n-1} 
	e ^{C \s_{n-1}^{-\m}}\,,
\end{equation}
for some constant $C>0$.

\item
For $n\ge 0$, there exists a sequence of functionals $F_n(u) \equiv F_n(\om, u(\omega))$, defined for 
$\om\in\calO^{(n-1)}$, of the form
\begin{equation}\label{effekappa}
F_n(u)=f_n + L_n u +Q_n(u),
\end{equation}
such that 

\begin{enumerate}

\item
$L_n$ is invertible for $\om\in\calO^{(n)}$ and setting
\begin{equation}\label{acchino}
h_n:=-L_n^{-1}f_n,
\end{equation}
there exists ${\mathtt r}_n ={\mathtt r}_{n}(\f)\in\calH (\T^\infty_{s_{n-1}-3\s_{n-1}})$ such that
\begin{equation}\label{riflettere}
\begin{aligned}
&F_{n}(u) = {\mathtt r}_nT_{n}^{-1}F_{n-1}(h_{n-1}+T_k u),\qquad n\ge1, \\
&|{\mathtt r}_n-1|_{s_{n-1}-3\s_{n-1}}^{\calO^{(n-1)}} \le \s_{n-1}^{-\tau_2} e^{C\s_{n-1}^{-\m}}\e_{n-1}
\end{aligned}
\end{equation}

\item
$f_n=f_n(\f,x)$ is a given function satisfying
\begin{equation}
|f_n|_{s_{n-1}-2\s_{n-1}}^{\calO^{(n-1)}}\lesssim \s_{n-1}^{-4} \e_{n-1}^2,\qquad n\ge1
\end{equation}

\item
 $L_n$ is a linear operator of the form
\begin{equation}\label{elleka}
L_n = \om\cdot \del_\f + (1+A_n)\del_{xxx} + B_n(\f,x)\del_x + C_n(\f,x)
\end{equation}
such that
\begin{equation}\label{media}
\frac{1}{2\pi}\int_\TTT B_n(\f,x) dx = \ol{b}_n
\end{equation}
and for $n\ge 1$
\begin{equation}\label{taglie}
\begin{aligned}
&|A_n-A_{n-1}|^{\calO^{(n-1)}} \le \s_{n-1}^{-\tau_2} e^{C\s_{n-1}^{-\mu}}\e_{n-1},\\
&|B_{n}-B_{n-1}|_{s_{n-1}-3\s_{n-1}}^{\calO^{(n-1)}}\lesssim \s_{n-1}^{-\tau_2}
e^{C\s_{n-1}^{- \mu}}\e_{n-1} \\
&|C_{n}-C_{n-1}|_{s_{n-1}-3\s_{n-1}}^{\calO^{(n-1)}}\lesssim \s_{n-1}^{-\tau_2}
e^{C\s_{n-1}^{- \mu}}\e_{n-1} \,.
\end{aligned}
\end{equation}

\item 
$Q_n$ is of the form
\begin{equation}\label{qk}
Q_n(u)= \sum_{0\le i\le 2,\,0\le j\le 3 \atop 0\le i+j\le 4} q^{(n)}_{i,j} (\f,x) (\del_x^i u)(\del_x^j u)
\end{equation}
with the coefficients $q^{(n)}_{i,j} (\f,x) $ satisfying \eqref{identica0} for $n=0$, while  for $n\ge 1$
\begin{equation}\label{identica}
\begin{aligned}
&\sum_{0\le i\le 2,\,0\le j\le 3 \atop 0\le i+j\le 4} |q^{(n)}_{i,j}|_{s_{n-1}-3\s_{n-1}}^{\calO^{(n-1)}}\le 
C\sum_{l=1}^{n}2^{-l}\,, \\
&|q^{(n)}_{i,j}-q^{(n-1)}_{i,j}|_{s_{n-1}-3\s_{n-1}}^{\calO^{(n-1)}}\lesssim \s_{n-1}^{-\tau_3}
e^{C\s_{n-1}^{- \mu }}\e_{n-1}\,.
\end{aligned}
\end{equation}

\end{enumerate}

\item
Finally one has
\begin{equation}\label{campane}
 |h_n|_{s_n}^{\calO^{(n)}}\le \e_n
\end{equation}

\end{enumerate}

Moreover, setting
\begin{equation}\label{infinito}
\calO^{(\io)}:=\bigcap_{n\ge0}\calO^{(n)},
\end{equation}
and
\begin{equation}\label{cheppa}
u_n = h_0 + \sum_{j=1}^n T_1 \circ\ldots\circ T_j h_j.
\end{equation}
then
$$
u_\io := \lim_{n\to\io} u_n
$$
is well defined for $\om\in\calO^{(\io)}$, belongs to 
$\calH_{\ol{s}}$,
and solves $F(u_\io)=0$.
Finally the $\calO^{(\io)}$ has positive measure; precisely
\begin{equation}\label{testa}
{\mathbb P}(\calO^{(\io)}) = 1 - O(\g_0)\,.
\end{equation}
\end{prop}

From Proposition \ref{ittero} our main result Theorem \ref{erteo} follows immediately by noting that 
\eqref{quantifico} and \eqref{piccolo0} follow from \eqref{certo} for an appropriate choice 
$\e(S-\ol{s})$.

\begin{rmk}\label{dellintro}
Let us spend few words on the strategy of the algorithm. At each step we apply an
affine change of variables translating the approximate solution to zero; the translation
is not particularly relevant and we perform it only to simplify the notation.
On the other hand the linear change of variables is crucial.

In \eqref{effekappa} we denote by $f_n$ the ``constant term'', by $L_n$ is the ``linearized''
term and by $Q_n$ the  ``quadratic'' part. In this way the approximate solution at the $n$-th
step is $h_n=-L_n^{-1}f_n$.

In a classical KAM algorithm, in order to invert $L_n$ one typically applies a linear
change of variables that diagonalizes $L_n$; this, together with the translation by $h_n$
is the affine change of variables mentioned above, at least in the classical KAM scheme.

Unfortunately, in the case of unbounded nonlinearities this cannot be done. Indeed in order
to diagonalize $L_n$ in the unbounded case, one needs it to be a pseudo-differential operator.
On the other hand, after the diagonalization is performed, one loses the
pseudo-differential structure for the subsequent step.
Thus we chose the operators $T_n$ in \eqref{tika} in such a way that we preserve the
PDE structure and at the same time we diagonalize the highest order terms.

In the \cite{BBM}-like algorithm the Authors do not apply any change of variables, but they
use the reducibility of $L_n$ only in order to deduce the estimates. However such a procedure
works only in Sobolev class. Indeed in the analytic case, at each iterative step one
needs to lose some analyticity, due to the small divisors. Since we are studying almost-periodic
solutions, we need the analytic setting to deal with the small divisors.
As usual, the problem is that
the loss of the analyticity is related to the size of the perturbation; in the present case, at each step
$L_n$ is a diagonal term plus a perturbation $O(\e_0)$ with the same $\e_0$ for all $n$.

A more refined approach is to consider $L_n$ as a small variation of $L_{n-1}$; however 
the problem is that such small variation is unbounded. As a consequence, the operators
$T_n$ are not ``close-to-identity''. However, since $F_n$ is a differential operator, then 
the effect of applying $T_n$ is simply a slight modification of the coefficients;
see \eqref{taglie} and \eqref{identica}. Hence there is a strong motivation for applying the operators
$T_n$. In principle we could have also diagonalized the terms up to order $-k$ for any $k\ge0$; however
the latter change of variables are close to the identity and they introduce pseudo-differential terms.
\end{rmk}

%%%%%%%%%%%%%%%%%%%%%%%%%%%%%%%%%%%%%%%%%%%%%%%%%%%%%
\subsection{The zero-th step}\label{zero}
Item $1., 2.$ are trivial for $n=0$ while item $3.(b),(c),(d)$ amount to the definition of $F_0$, see 
\eqref{Fzero},\eqref{tayzero},\eqref{q0}.
Regarding item $3.(a)$ the invertibility of $L_0$ follows from the definition of $\calO^{(0)}$.
Indeed, consider the equation
\begin{equation}\label{h0}
L_0 h_0 = -f_0
\end{equation}
with
\[
\avg{f_0(\f,\cdot)}_{x} =0
\]
we have the following result.

\begin{lemma}[\bf Homological equation] \label{stac}
	Let $s > 0$, $0<\s<1$, $f_0 \in\calH_{s+\s}$, $\omega \in\calO^{(0)}$ (see \eqref{diofantinoBIS}).
	Then there exists a unique solution 
	$h_0   \in \calH_s$ of \eqref{h0}	.
	Moreover one has
	$$
	| h_0 |_s^{\calO^{(0)}} \lesssim \g^{-1}{\rm exp}\Big(\frac{\tau}{\sigma^{\frac{1}{\zia}}} 
	\ln\Big(\frac{\tau}{\sigma} \Big) \Big) 
	| f |_{s + \sigma }\,. 
	$$
	for some constant $\tau=\tau(\zia) > 0$. 
\end{lemma}

\begin{rmk}\label{vini}
Note that from Lemma \ref{stac} above it follows that there is $\tC_0$ such that
a solution $h_0$ of \eqref{h0} actually satisfies
\begin{equation}\label{vento}
	| h_0 |_s^{\calO^{(0)}} \lesssim e^{\tC_0 {\sigma^{-\mu}} } 
	| f |_{s + \sigma }\,. 
\end{equation}
where we recall that by \eqref{costanti}, $\mu > \max\{ 1 , \frac{1}{\eta}\}$. Of course the constant $\tC_0$ 
is correlated with the correction to the exponent $\frac{1}{\h}$.
\end{rmk}

From Lemma \ref{stac} and \eqref{h0} it follows that
$h_0$ is analytic in a strip $s_0$ (where $S=s_0+\s_{-1}$ is the analyticity of $f$, to be chosen).
Moreover, by Lemma \ref{stac}
the size of $h_0$ is
\begin{equation}\label{elabase}
|h_0|_{s_0}^{\calO^{(0)}}\sim  e^{\tC_0\sigma_{-1}^{-\m}} |f_0|_{S}
\end{equation}
proving item 4. for $|f_0|_S$ small enough, which is true by \eqref{quantifico}.

\subsection{The $ n+1$-th step}\label{esimo}

Assume now that we iterated the procedure above up to $n\ge 0$ times. 
This means that we arrived at a quadratic equation
\begin{equation}\label{ekappa}
F_n(u)=0,\qquad
F_n(u) = f_n + L_n u + Q_n(u).
\end{equation}
Defined on $\calO^{(n-1)}$ (recall that $\calO^{(-1)}=\Dc$).

By the inductive hypothesis \eqref{identica} we deduce
 that for all $0<s-\s < s_{n-1}- 3\s_{n-1}$ one has
\begin{subequations}\label{qkappa}
\begin{align}
&|Q_n(u)|_{s-\s}^{\calO^{(n-1)}} \lesssim  \s^{-4}(|u|_s^{\calO^{(n-1)}})^2 \label{q1k} \\
&|Q_n'(u)[h]|_{s-\s}^{\calO^{(n-1)}}\lesssim \s^{-4}|u|_s^{\calO^{(n-1)}} |h|_s^{\calO^{(n-1)}}\label{q2k}
\end{align}
\end{subequations}

Moreover, again by the inductive hypothesis,  we can invert $L_n$ and define $h_n$ by \eqref{acchino}.
Now we set
\begin{equation}\label{corona}
F_{n+1}(v) = {\mathtt r}_{n+1} T_{n+1}^{-1}F_n(h_n+ T_{n+1}v )
\end{equation}
where 
\begin{equation}\label{nonsicak}
T_{n+1} v(\f,x)=  (1+\xi^{(n+1)}_x)v( \f +\omega \be^{(n+1)}(\f), x+ \xi^{(n+1)}(\f,x) + p^{(n+1)}(\f))\,
\end{equation}
and $r_{n+1}$ are to be chosen in order to ensure that $L_{n+1}:= F'_{n+1}(0)$ has the form \eqref{elleka} 
with $n\rightsquigarrow n+1$.

Of course by Taylor expansion we can identify
\begin{equation}\label{stronzi}
\begin{aligned}
f_{n+1}&=  {\mathtt r}_{n+1}T_{n+1}^{-1}(f_n + L_n(h_n) + Q_n(h_n))= 
 {\mathtt r}_{n+1}T_{n+1}^{-1} Q_n(h_n)\,,\\ 
L_{n+1} &=  {\mathtt r}_{n+1}T_{n+1}^{-1}(L_n + Q'_n(h_n)) T_{n+1}\,\\
Q_{n+1}(v) &=  {\mathtt r}_{n+1}(T_{n+1}^{-1} (Q_n( h_n + T_{n+1} v) - Q_n(h_n)- Q'_n(h_n) T_{n+1} v)) \\
&=  {\mathtt r}_{n+1}T_{n+1}^{-1}Q_n(T_{n+1} v)\,.
\end{aligned}
\end{equation}

\begin{rmk}
Note that the last equality in \eqref{stronzi} follows from the fact that the nonlinearity $Q$ in \eqref{kdv}
is  quadratic. In the general case, the last term is controlled by the second derivative, and thus 
one has to assume a bound of the type \eqref{qkappa} for $Q''$.
\end{rmk}

In section \ref{dimo} we prove the following
\begin{prop}\label{comecaz}
	Assuming that
	\begin{equation}
	\label{smallk}
	\e_n \le \sigma_n^{\tau_1 +1}e^{-C\sigma_n^{-\mu}}
	\end{equation}
	for some $C>0$,
	there exist $\xi^{(n+1)}$, $\beta^{(n+1)}$, $p^{(n+1)}$ and $ {\mathtt r}_{n+1}\in\calH(\T^\infty_{s_n - \s_n} \times \T_{s_n - \s_n})$, 
	defined for all $\omega\in {\calO^{(n)}}$  and satisfying
	\begin{equation}\label{tagliere}
	| \xi^{(n+1)}|^{\calO^{(n)}}_{s_n- \s_n},|\be^{(n+1)}|^{\calO^{(n)}}_{s_n-\s_n},
	|p^{(n+1)}|^{\calO^{(n)}}_{s_n-\s_n} ,
	|{\mathtt r}_{n+1}-1|^{\calO^{(n)}}_{s_n-{\s_n}} 
	\lesssim \s_n^{-\tau_1} \e_n e ^{C \s_n^{- \mu}}\,
	\end{equation}
	such that \eqref{nonsicak} is well defined and symplectic as well as its inverse, and moreover
	\begin{equation}\label{equello}
	{\mathtt r}_{n+1}T_{n+1}^{-1}(L_n + Q'_{n}(h_n))T_{n+1} = \om\cdot \del_\f + (1+A_{n+1})\del_{xxx} + 
	B_{n+1}(\f,x)\del_x + C_{n+1}(\f,x)
	\end{equation}
	and \eqref{media} and \eqref{taglie} hold with $n\rightsquigarrow n+1$.
\end{prop}

The assumption \eqref{smallk} follows from \eqref{piccolo0}, provided that we choose 
the constants $\tau,\tC$ and $\epsilon_0$ appropriately. %Indeed \textr{BLABLABLA}

We now prove \eqref{qk} and \eqref{identica} for $n\rightsquigarrow n+1$, namely the following
result.

\begin{lemma}\label{pranzo}
	One has
	\begin{equation}\label{qkn}
	Q_{n+1}(v)={\mathtt r}_{n+1} T_{n+1}^{-1}Q_n(T_{n+1} v)= {\mathtt r}_{n+1}
	\sum_{0\le i\le 2,\,0\le j\le 3 \atop 0\le i+j\le 4} q^{(n+1)}_{i,j} (\f,x) (\del_x^i v)(\del_x^j v)
	\end{equation}
	with the coefficients $q^{(n+1)}_{i,j} (\f,x) $ satisfying
	\begin{equation}\label{identica+}
	\begin{aligned}
	&\sum_{0\le i\le 2,\,0\le j\le 3 \atop 0\le i+j\le 4} |q^{(n+1)}_{i,j}|^{\calO^{(n)}}_{s_{n}-3\s_{n}}\le 
	C\sum_{l=1}^{n+1}2^{-l}\,, \\
	&|q^{(n+1)}_{i,j}-q^{(n)}_{i,j}|^{\calO^{(n)}}_{s_{n}-3\s_{n}}\lesssim
	 \s_{n}^{-\tau_3}e^{C\s_{n}^{-\mu}}\e_{n}\,.
	\end{aligned}
	\end{equation}
\end{lemma}

\prova
By construction
\begin{equation}\label{colti}
Q_{n+1}(u)={\mathtt r}_{n+1} \sum_{0\le i\le 2,\,0\le j\le 3 \atop 0\le i+j\le 4} 
T_{n+1}^{-1}[q^{(n)}_{i,j} (\f,x) (\del_x^i T_{n+1}v)(\del_x^j T_{n+1} v)].
\end{equation}
Now we first note that
\[
\del_x (T_{n+1}v) = \x^{(n+1)}_{xx} v(\theta,y) + (1+\x_x)^2v_{y}(\theta,y)
\]
where
$$
(\theta,y)=(\f +\omega \be^{(n+1)}(\f), x+ \xi^{(n+1)}(\f,x) + p^{(n+1)}(\f)).
$$
Hence
the terms $\del_x^i T_{n+1}v$ are of the form
\begin{equation}\label{abbio}
\del_x^i T_{n+1}v = \del_y^i v(\theta,y) + \sum_{l=0}^{i} g_{l,i}(\f,x)\del_y^l v(\theta,y),\qquad
%g_{l,i}(\f,x) = \left\{
%\begin{aligned}
%&1+ O(\x),\quad l=i \\
%&O(\x),\quad l<i.
%\end{aligned}
%\right.
|g_{l,i}|_{s_n-2\s_n}^{\calO^{(n)}} \lesssim \s_n^{-(i+2)}|\x^{(n+1)}|_{s_n-\s_n}^{\calO^{(n)}}
\end{equation}
%where by $O(\x)$ we mean that it is a function depending on $\x$ and its derivatives up to order $i+1$. Inserting
%Note that by the smallness hypothesis $|1+\x_x^{(k+1)}|<2$.
\\
Inserting
\eqref{abbio} into \eqref{colti} we get 
\begin{equation}\label{cco}
\begin{aligned}
q^{(n+1)}_{l,m}&= {\mathtt r}_{n+1}\left( T_{n+1}^{-1} q^{(n)}_{l,m} +
\sum_{j=0}^4T_{n+1}^{-1}(q_{l,j}^{(n)}g_{m,j}) + \sum_{i=0}^4T_{n+1}^{-1}(q_{i,m}^{(n)}g_{l,i})\right.\\
&\qquad+\left.
\sum_{0\le i\le 2,\,0\le j\le 3 \atop 0\le i+j\le 4} T_{n+1}^{-1}(q^{(n)}_{i,j}g_{l,i}g_{m,j}) 
\right)%(\del_y^l v)( \del_y^m v).
\end{aligned}
\end{equation}
so that
\begin{equation}\label{sonno}
q^{(n+1)}_{i,j}= T_{n+1}^{-1}(q^{(n)}_{i,j} + O(\x_{n+1})),\qquad 
|T^{-1}_{n+1}O(\x_{n+1})|_{s_{n}-3\s_n}^{\calO^{(n)}}
\lesssim\s_n^{-\tau_3}\e_ne^{C\s_n^{\mu}}.
\end{equation}
In order to obtain the bound \eqref{sonno} we used the first line of \eqref{identica}
 to control the sums appearing in \eqref{cco}.

Finally , since
\[
T_{n+1}^{-1}(q)- q:= (1+\tilde\xi^{(n+1)}_x)q( \f,x)- q(\theta,y)
\]
the bound follows.
\EP

Now,
by  \eqref{q1k} and \eqref{stronzi}  $f_{n+1}=f_{n+1}(\f,x)$ satisfies
\begin{equation}\label{ronf}
|f_{n+1}|_{s_{n}-2\s_{n}}^{\calO^{(n)}}\lesssim \s_{n}^{-4} \e_{n}^2.
\end{equation}

In Section \ref{eppoi} we prove the existence of a Cantor set $\calO^{(n+1)}$ where item 
$3.(a)$ of the iterative lemma  
holds with $n\rightsquigarrow n+1$.

\begin{prop}\label{ciaraggio}

Assume  that 
\begin{equation}\label{piccino}
2^{n}\s_n^{-\tau}e^{\tC \s_n^{-\mu}}\e_n\ll1 \,,
\end{equation}
with $\tau\ge\tau_2$.
Setting $\la_3^{(n+1)}:=1+A_{n+1}$,
	there exist Lipschitz functions 
	\begin{equation}\label{omegone+}
	\Omega^{(n+1)} (j)= \lambda^{(n+1)}_3 j^3 + \lambda^{(n+1)}_1 j + r^{(n+1)}_j
	\end{equation}
	satisfying
	\begin{equation}\label{tutto+}
	|\la_1^{(n+1)}-\la_1^{(n)}|^{\calO^{(n)}},
	\sup_{j\in\ZZZ\setminus\{0\}}|r_j^{(n+1)}-r_j^{(n)}|^{\calO^{(n)}}\lesssim 
	\s_{n}^{-\tau}\e_{n}e^{C\s_{n}^{- \mu}}
	\end{equation}
	such that setting
	\begin{equation}\label{caloni+}
	\calE^{(n+1)}:=\Big\{\om\in\calO^{(n)}:\; |\om\cdot\ell + \Omega^{(n+1)}(j)-\Omega^{(n+1)}(h)|\ge 
	\frac{2\gamma_{n+1}|j^3 - h^3|}{{\mathtt d}(\ell)}\,,
		\ \forall %\ell\in \Z^\io_*\,,\;  h,j\in \N \,,\; 
		(\ell,h,j)\ne (0,h,h) \Big\}
	\end{equation}
	for $\om\in\calE^{(n+1)}$ there exists an invertible and bounded linear operator $M^{(n+1)}$
	\begin{equation}\label{ukappa+}
	\|M^{(n+1)}-\Id \|_{s_{n}-5\s_{n}}^{\calE^{(n+1)}}\le \s_{0}^{-\tau}e^{\tC \s_{0}^{-\mu}}\e_{0}\\
	\end{equation}
	such that 
	\begin{equation}\label{diago+}
	(M^{(n+1)})^{-1} L_{n+1} M^{(n+1)} = D_{n+1}= 
	{\rm diag}\pa{\omega\cdot \ell + \Omega^{(n+1)}(j) }_{(\ell,j)\in \Z^\io_*\times\ZZZ\setminus\{0\}}
	\end{equation}
	
\end{prop}

The assumption \eqref{piccino} follows from \eqref{piccolo0}, provided that we choose 
the constants $\tau,\tC$ and $\epsilon_0$ appropriately. %Indeed \textr{BLABLABLA}

\begin{rmk}\label{lemn}
Note that in the context of \cite{CFP}
Proposition \ref{ciaraggio} is much simpler to prove, because in order to diagonalize
the linearized operator one uses tame estimates coming from the Sobolev regularity on
the boundary of the domain. Then the smallness conditions are much simpler to handle.
Here we have to strongly rely on the fact that $L_{n+1}$ is a ``small'' unbounded perturbation
of $L_n$ in order to show that
the operators $M^{(n)}$ and $M^{(n+1)}$ are close to 
each other.  This is a very delicate
issue; see Lemma \ref{lemma variazione parte partial x} and Section \ref{meloinculo},
which are probably the more technical parts of this paper. 
\end{rmk}

\begin{lemma}[\bf Homological equation] \label{eppa}
Set
\begin{equation}\label{mel1}
	\calU^{(n+1)}:=\set{\om\in\calO^{(n)}:\quad |\om\cdot\ell + \Omega^{(n+1)}(j)|\ge \gamma_{n+1}
		\frac{|j|^3}{{\mathtt d}(\ell)}\,,
		\quad \forall %\ell\in \Z^\io_*\,,\quad  j\in \N \,,\; 
		(\ell,j)\ne (0,0) }
\end{equation}
	For $\omega \in\calO^{(n+1)}:= \calU^{(n+1)}\cap\calE^{(n+1)}$ one has
	\begin{equation}\label{attimo}
	h_{n+1}:= - L_{n+1}^{-1}f_{n+1} \in  \calH_{s_{n+1}}
	\end{equation}
	and one has
	$$
	| h_{n+1} |_{s_{n+1}}^{\calO^{(n+1)}} \lesssim {\rm exp}
	\Big({\tau}{\sigma_{n}^{-\frac{1}{\zia}}} \ln\Big(\frac{\tau}{\sigma_n} \Big) \Big) 
	| f_{n+1} |_{s_{n+1} + \sigma_n }^{\calO^{(n)}}\,. 
	$$
\end{lemma}

\proof
The result follows simply by using the definition of $\calO^{(n+1)}$ and applying 
Lemma \ref{cazzodilemma}.
\EP

Of course from Lemma \ref{eppa} it follows that,
\begin{equation}\label{pipino}
|h_{n+1}|_{s_{n+1}}^{\calO^{(n+1)}} \lesssim \s_{n}^{-4}e^{C\s_{n}^{- \mu}}\e_{n}^2
\end{equation}

Now we want to show inductively that
\begin{equation}\label{super}
\s_{n}^{-4}e^{C\s_{n}^{- \mu}}\e_{n}^2
 \le \e_0 e^{-\chi^n+1},\qquad \qquad \chi=\frac{3}{2}
\end{equation}
for $\e_0$ small enough.

By the definition of $\e_n$ in \eqref{costanti}, \eqref{super} is equivalent to
%\begin{equation}\label{aspe}
% \s_{k}^{-4}e^{C\s_{k}^{-\frac{1}{\h}+}}  \e_0^2 e^{-2\chi^k+2} \lesssim \e_0 e^{-\chi^{k+1}+1} 
%\end{equation}
%which in turn is equivalent to
\begin{equation}\label{sexy}
\e_0  \lesssim \s_0^4 n^{-8}  e^{ \chi^{n}(2-\chi)-C'{n^{\mu}}}
\end{equation}

Since the r.h.s. of \eqref{sexy} admits a positive minimum, we can regard it as a smallness condition on $\e_0$,
which is precisely \eqref{piccolo0}.

We now prove \eqref{alvolo} with $n\rightsquigarrow n+1$.
We only prove the bound for the set ${\cal E}^{(n)} \setminus {\cal E}^{(n+1)}$. 
The other one can be proved by similar arguments (it is actually even easier). 
Let us start by  writing

\begin{equation}\label{def insiemi risonanti}
\begin{aligned}
&{\cal E}^{(n)} \setminus {\cal E}^{(n+1)} = \bigcup_{(\ell, j, j') \neq (0, j, j)} {\cal R}(\ell, j , j')\,, \\
& {\cal R}(\ell, j , j') := \Big\{ \omega \in {\cal E}^{(n)} : 
|\omega \cdot \ell + \Omega^{(n+1)}(j) - \Omega^{(n+1)}(j')| < 
\frac{2 \gamma_{n+1} |j^3 - j'^3|}{\mathtt d(\ell)} \Big\}\,, \\
& \forall (\ell, j, j' ) \in \Z^\infty_* \times (\Z  \setminus \{ 0 \})\times (\Z  \setminus \{ 0 \}), \quad (\ell, j, j') \neq (0, j, j)\,. 
\end{aligned}
\end{equation}

\begin{lemma}\label{lemma risonanti 1}
Denote $|\ell|_1$ as in \eqref{Z inf *} with $\h\rightsquigarrow 1$.
 For any $(\ell, j, j') \neq (0, j, j)$ such that $|\ell|_1 \leq n^{2}$, one has that ${\cal R}(\ell, j , j') = \emptyset$.
\end{lemma}

\prova
Let $(\ell, j, j') \in \Z^\infty_* \times (\Z  \setminus \{ 0 \})\times (\Z  \setminus \{ 0 \})$, 
$(\ell, j, j') \neq (0, j, j)$, $|\ell|_1 \leq n^{2}$. If $j = j'$, clearly $\ell \neq 0$ and 
${\cal R}(\ell, 0, 0) = \emptyset$ because $\omega\in\Dcb$ with $\ol{\g}>2 \gamma_{n+1}$; recall 
\eqref{costanti}.
Hence we are left to 
analyze the case  $j \neq j'$. 

By  \eqref{tutto+},  for any $j, j' \in \Z \setminus \{ 0 \}$, $j \neq j'$
\begin{equation}\label{blo corona}
\Big| \Big( \Omega^{(n+1)}(j) - \Omega^{(n+1)}(j') \Big) - 
\Big(\Omega^{(n)}(j) - \Omega^{(n)}(j') \Big) \Big|\\
%& \leq  |\lambda_3^+ - \lambda_3| |j^3 - j'^3|+ |\lambda_1^+ - \lambda_1||j - j'| + |q_n^+(j)| +  |q_n(j)| \\
 {\lesssim}  \,
\s_{n}^{-\tau}\e_{n}e^{C\s_{n}^{- \mu}} |j^3 - j'^3|\,.
\end{equation}

Therefore, for any $\omega \in {\cal E}^{(n)}$
\begin{equation}\label{blo corona 2}
\begin{aligned}
|\omega \cdot \ell + \Omega^{(n+1)}(j) - \Omega^{(n+1)}(j')| & 
\geq |\omega \cdot \ell + \Omega^{(n)}(j) - \Omega^{(n)}(j')| \\
&\qquad\qquad- 
\Big| \Big( \Omega^{(n+1)}(j) - \Omega^{(n+1)}(j') \Big) - 
\Big(\Omega^{(n)}(j) - \Omega^{(n)}(j') \Big) \Big|  \\
& \geq \frac{2 \gamma_n |j^3 - j'^3|}{\mathtt d(\ell)} - C 
 \s_{n}^{-\tau}\e_{n}e^{C\s_{n}^{- \mu}} |j^3 - j'^3| \\
 &\geq \frac{2 \gamma_{n+1}|j^3 - j'^3|}{\mathtt d(\ell)}
\end{aligned}
\end{equation} 
where in the last inequality we used \eqref{costanti} and the fact that, by \eqref{taglio} one has
\[
\begin{aligned}
\s_{n}^{-\tau}\e_{n}e^{C\s_{n}^{- \mu}}{\mathtt d(\ell)} &\le
\s_{n}^{-\tau}\e_{n}e^{C\s_{n}^{- \mu}} (1+n^2)^{C(1)n}\le \g_0 2^{-n}\,.
\end{aligned}
\]

The estimate \eqref{blo corona 2} clearly implies that ${\cal R}(\ell, j, j') = \emptyset$ for $|\ell|_1 \leq n^2$. 
\EP

\begin{lemma}\label{lemma risonanti 2}
Let ${\cal R}(\ell, j, j') \neq \emptyset$. Then $\ell \neq 0$, $|j^3 - j'^3|
 \lesssim \| \ell \|_1$ and ${\mathbb P}\Big({\cal R}(\ell, j, j')\Big) \lesssim \frac{\gamma_{n+1}}{\mathtt d(\ell)}$
\end{lemma}

\prova
The proof is identical to the one for Lemma 6.2 in \cite{MP}, simply replacing $j^2$ with $j^3$.
\EP

By \eqref{def insiemi risonanti} and collecting Lemmata 
\ref{lemma risonanti 1}, \ref{lemma risonanti 2}, one obtains that
\begin{equation}
\begin{aligned}
{\mathbb P}\Big( {\cal E}^{(n)} \setminus {\cal E}^{(n+1)}\Big) & \lesssim \sum_{\begin{subarray}{c}
|\ell|_1 \geq n^2 \\
|j|, |j'| \leq C \| \ell \|_1 
\end{subarray}}  \frac{\g_{n+1}}{\mathtt d(\ell)} \lesssim\g_{n+1} 
\sum_{|\ell|_1 \geq n^2} \frac{\| \ell \|_1^2}{\mathtt d(\ell)} \lesssim \g_{n+1} n^{- 2} 
\sum_{\ell \in \Z^\infty_*} \frac{|\ell|_1^3}{\mathtt d(\ell)}  \lesssim \g_{n+1} n^{- 2}.
\end{aligned}
\end{equation} 
where in the last inequality we used Lemma \ref{palline}.
% that the series
%$$
%\sum_{\ell \in \Z^\infty_*} \frac{|\ell|_1^3}{\mathtt d(\ell)}
%$$
%is convergent; see .
Thus \eqref{alvolo} follows.
%{\color{red} DOVREBBE ESSERE SIMILE AL LEMMA  B.3 in \cite{MP}}

We now the convergence of the scheme.
Precisely we show that the series \eqref{cheppa} converges totally in $\calH_{\ol{s}}$ .
Note that 
\begin{equation}\label{hocca}
|T_i u|_{\ol s}^{\calO^{(\io)}} \le (1+2^{-i})|u|_{\ol{s}+\s_i}^{\calO^{(\io)}} \le 2|u|_{\ol{s}+\s_i}^{\calO^{(\io)}}.
\end{equation}
Thuse, using \eqref{hocca} into \eqref{cheppa} we get
\begin{equation}\label{esce}
|u_n|_{\ol{s}}^{\calO^{(\io)}} \le |h_0|_{\ol{s}}^{\calO^{(\io)}} + \sum_{j=1}^n 2^j 
|h_j|_{\ol{s}+ (\s_1+\ldots +\s_j)}^{\calO^{(\io)}}
\end{equation}
Now since
\begin{equation}\label{vivi}
\ol{s} + \sum_{n=1}^{\io}\s_n = \ol{s} + \frac{6\s_{-1}}{\pi^2} \sum_{n\ge1} \frac{1}{n^2} = s_\io \le s_j
\end{equation}
we deduce that $u_\io\in\calH_{\ol{s}}$. Finally by continuity
$$
F(u_\io) = \lim_{n\to\io} F(u_n) =\lim_{n\to\io} T_1^{-1}T_2^{-1}\ldots T_n^{-1} F_n(h_n)=0.
$$
so the assertion follows since (recall $\ol{s}:= s_\io - \sum_{n\ge1} \s_n$  and \eqref{vivi})
\[
|T_1^{-1}T_2^{-1}\ldots T_n^{-1} F_n(h_n)|_{\ol{s}}^{\calO^{(\io)}} \le 2^n \s_n^{-4}\e_n^2 \,.
\]

We finally conclude the proof of Proposition \ref{ittero} by showing that \eqref{testa} holds.

First of all, reasoning as in Lemma \ref{lemma risonanti 2} and using 
Lemma \ref{palline}, we see that 
\[
{\mathbb P}(\calO^{(0)}) = 1-O(\g_0)
\]
Then 
\[
{\mathbb P}(\calO^{(\io)}) =  {\mathbb P}(\calO^{(0)}) - \sum_{n\ge0}{\mathbb P}(\calO^{(n)}
\setminus\calO^{(n+1)})
\]
so that \eqref{testa} follows by \eqref{alvolo}.
\EP

%%%%%%%%%%%%%%%%%%%%%%%%%%%%%%%%%%%%%%%%%%%%%%%%%
%%%%%%%%%%%%%%%%%%%%%%%%%%%%%%%%%%%%%%%%%%%%%%%%%
\section{Proof of Proposition \ref{comecaz}}\label{dimo}
\zerarcounters
%%%%%%%%%%%%%%%%%%%%%%%%%%%%%%%%%%%%%%%%%%%%%%%%%
%%%%%%%%%%%%%%%%%%%%%%%%%%%%%%%%%%%%%%%%%%%%%%%%%

In order to prove Proposition \ref{comecaz}, we start by dropping the index $n$, i.e. we set
${\cal L} \equiv L_n$ 
(see \eqref{elleka}) and ${\cal Q} \equiv Q'_n(h_n)$ (see \eqref{stronzi}).

More generally, we consider a Hamiltonian 
operator of the form 
\begin{equation}\label{operatore cal L}
\begin{aligned}
& {\cal L}^{(0)} = {\cal L} + {\cal Q}  \\
& {\cal L}: =  \omega \cdot \partial_\f + \lambda_3 \partial_x^3 + a_1(\f, x ) \partial_x + a_0(\f, x)\,, \\
& {\cal Q} :=  d_3(\f, x) \partial_x^3 + d_2(\f , x) \partial_x^2 + d_1(\f,  x) \partial_x + d_0(\f , x)
\end{aligned}
\end{equation}
defined for all $\om\in\Omega\subseteq \Dc$ and
 $\lambda_3$, $a_0, a_1, d_0 \ldots, d_3$ satisfy the following properties.

\begin{enumerate}

\item There is $\de_0$ small enough such that
{\begin{equation}\label{stime lambda 3Mic}
	|\lambda_3 -1|^{\Omega} \le \delta_0
\end{equation}
}

\item There is $\rho>0$ such that $a_i\in\calH(\T^\infty_\rho\times \T_\rho)$ and
{\begin{equation}\label{vaffanculo}
 | a_i |_{\rho}^{\Omega} \le  \delta_0\,,\qquad i=0,1
\end{equation}
and moreover
{\begin{equation}\label{prop mediaMic}
 \lambda_1 := \frac{1}{2 \pi} \int_\T a_1(\f, x)\, d x 
\end{equation}}
i.e. it does not depend on $\f$.
}

\item
$d_0 \ldots, d_3 \in \calH(\T^\infty_\rho\times \T_\rho)$ 
(note that by the Hamiltonian structure $d_2 = \partial_x d_3$) and they satisfy the estimate 
\begin{equation}\label{stima a3Mic}
| d_i |_{\rho}^{\Omega}  \lesssim \delta\,,
\end{equation}  
for some $\de \ll \min\{\de_0,\rho\}$.

\end{enumerate} 

Let us now choose $\ze$ such that $0<\ze\ll\rho$ and
\begin{equation}\label{ruota}
\zeta^{- \tau'} e^{2 \tC_0 \zeta^{- \mu }} \delta\ll1.
\end{equation}
for some $\tau'>0$. We shall conjugate $\LL^{(0)}$ to a new operator $\frac{1}{\mathtt r}\LL_+$
with ${\mathtt r}={\mathtt r}(\f)$ an explicit function with 
\begin{equation}\label{elle+}
\LL_+ = \om\cdot\del_\f +\la_3^+\del_x^3 + a_1^+(\f,x)\del_x + a_0^+(\f,x)
\end{equation}
with the coefficients satisfying 
\begin{equation}\label{ufo}
|\lambda_3^+ -\la_3|^{\Omega} \lesssim \de
\end{equation}
and
\begin{equation}\label{robot}
 | a_i^+ - a_i |_{\rho-2\ze}^{\Omega} \le  \zeta^{- \tau'} e^{2 \tC_0 \zeta^{- \mu }} \delta
 ,\qquad
  \lambda_1 := \frac{1}{2 \pi} \int_\T a_1(\f, x)\, d x .
\end{equation}

This will allow us to conclude the proof of Proposition \ref{comecaz}.

%%%%%%%%%%%%%%%%%%%%%%%%%%%%%%%%%%%%%%%%%%%%%%%%%
\subsection{Elimination of the $x$-dependence from the highest order term}\label{gatto}
%%%%%%%%%%%%%%%%%%%%%%%%%%%%%%%%%%%%%%%%%%%%%%%%%

Consider an analytic function $\alpha (\f, x)$ (to be determined) and let 
$$
{\cal T}_1 u (\f, x) := (1 + \alpha_x(\f, x)) ({\cal A} u)(\f, x), \quad {\cal A} u(\f, x) :=  u(\f, x + \alpha(\f,  x))\,. 
$$
We choose $\alpha(\f, x)$ and $m_3(\f)$ in such a way that
\begin{equation}\label{eq omologica grado 3}
(\lambda_3 + d_3(\f, x)) \big( 1 + \alpha_x(\f, x)\big)^{3} =  m_3(\f)\,,
\end{equation} 
which implies
\begin{equation}\label{scelta beta}
\alpha(\f, x) := \partial_x^{- 1} \Big[ \frac{m_3(\f)^{\frac13}}{\big(\lambda_3 + d_3(\f, x) \big)^{\frac13}} - 1 \Big], 
\quad m_3(\f) := \Big(\frac{1}{2 \pi} \int_\T \frac{d x}{\big( \lambda_3 + d_3(\f, x) \big)^{\frac13}} \Big)^{- 3}\,.
\end{equation}
By \eqref{stime lambda 3Mic}, \eqref{stima a3Mic} and Lemma \ref{moser type lemma} one has 
\begin{equation}\label{stime m3 alpha}
| m_3 - \lambda_3 |_{\rho}^{\Omega}\,,\, | \alpha |_{\rho}^{\Omega} \lesssim \delta 
\end{equation}

Note that for any  $0<\zeta\ll \rho$ such that  $\delta \zeta^{-1} \ll 1$, by Lemma 
\ref{lemma diffeo inverso}, $x \mapsto x + \alpha(\f, x)$ is invertible and the inverse 
is given by $y \mapsto y + \widetilde \alpha(\f, y)$ with 
\begin{equation}\label{bellaciao}
\widetilde \alpha \in {\cal H}(\T^\infty_{\rho - \zeta} \times \T_{\rho - \zeta}), \quad
 | \widetilde \alpha |_{\rho - \zeta }^{\Omega} , | \alpha |_{\rho}^{\Omega} \lesssim \delta\,. 
\end{equation}

A direct calculations shows that 
\begin{equation}\label{inverse trnaspose of cal A}
\begin{aligned}
& {\cal A}^{- 1} u(\f, y) =  u(\f, y + \widetilde \alpha(\f, y))\,,  
 \quad  {\cal T}_1^{- 1} = (1 + \widetilde \alpha_y) {\cal A}^{- 1}
\end{aligned}
\end{equation}
and the following conjugation rules hold:
\begin{equation}\label{coniugazione cal A}
\begin{aligned}
& {\cal T}_1^{- 1}  \, a(\f, x) \,  {\cal T}_1  = {\cal A}^{- 1}  \, a(\f, x) \,  {\cal A} =  ({\cal A}^{- 1} a)(\f, y)\,, \\
& {\cal T}_1^{- 1} \partial_x {\cal T}_1 = \Big( 1 + {\cal A}^{- 1}(\alpha_x) \Big) \partial_y +  
 (1 + \widetilde \alpha_y) {\cal A}^{- 1}(\alpha_{xx}) \,, \\
& {\cal T}_1^{- 1} \omega \cdot \partial_\f {\cal T}_1 = \omega \cdot \partial_\f + 
{\cal A}^{- 1}(\omega \cdot \partial_\f \alpha) \partial_y + (1 + \widetilde \alpha_y)
 {\cal A}^{- 1}(\omega \cdot \partial_\f \al_x)\,. 
\end{aligned}
\end{equation}

Clearly one can get similar conjugation formulae for higher order derivatives,
having expression similar to \eqref{abbio}.
In conclusion
\begin{equation}\label{cal L (1)}
\begin{aligned}
 {\cal L}^{(1)} &:= {\cal T}_1^{- 1} ({\cal L} + {\cal Q}) {\cal T}_1 \\
&= \omega \cdot \partial_\f + {\cal A}^{- 1}\Big[ (\lambda_3 + q_3)(1 + \al_x)^3 \Big] \partial_y^3 
+ b_2(\f, y) \partial_y^2 + b_1(\f, y) \partial_y + b_0(\f, y)\\
&=\omega \cdot \partial_\f + m_3(\f) \partial_x^3 + b_1(\f, x) \partial_x + b_0(\f, x)
\end{aligned}
\end{equation} 
for some (explicitly computable) coefficients $b_i$,
where in the last equality we used \eqref{eq omologica grado 3}
and the fact that ${\cal T}_1$ is symplectic, so that 
$b_2(\f, x) = 2 \partial_x m_3(\f) = 0$.

Furthermore, the estimates \eqref{stime lambda 3Mic}, \eqref{vaffanculo},
 \eqref{stime m3 alpha}, \eqref{bellaciao}, Corollary \ref{usiamoquesto}
 and Lemmata \ref{algebra analitiche}, 
 \ref{stime cauchy} imply that for $0 < \zeta \ll \rho$
\begin{equation}\label{bi - ai}
| b_i |_{\rho - 2\zeta}^{\Omega} \lesssim \de_0, \quad | b_i - a_i |_{\rho -2\zeta}^{\Omega}
 \lesssim \zeta^{- \tau} \delta, 
\quad \text{for some} \quad \tau > 0\,. 
\end{equation}

%%%%%%%%%%%%%%%%%%%%%%%%%%%%%%%%%%%%%%%%%%%%%%%%%
\subsection{Elimination of the $\f$-dependence from the highest order term}
%%%%%%%%%%%%%%%%%%%%%%%%%%%%%%%%%%%%%%%%%%%%%%%%%

We now consider a quasi periodic reparametrization of time of the form 
\begin{equation}\label{df cal B}
{\cal T}_2 u(\f, x) := u(\f + \omega \beta(\f), x)
\end{equation}
where $\beta : \T^\infty_{\rho-\ze} \to \R$ is an analytic function to be determined. Precisely we 
choose $\lambda_3^+ \in \R$ and $\beta(\f)$ in such a way that 
\begin{equation}\label{eq omologica riparametrizzazione tempo}
\lambda_3^+ \Big(1 + \omega \cdot \partial_\f \beta(\f) \Big) = m_3(\f)\,,
\end{equation}
obtaining thus
\begin{equation}\label{def }
\lambda_3^+ := \int_{\T^\infty} m_3(\f)\,d \f, \quad \beta(\f) := 
(\omega \cdot \partial_\f)^{- 1}\Big[ \frac{m_3}{\lambda_3^+} - 1\Big] 
\end{equation}
where we recall the definition \ref{media toro infinito dimensionale}. By the estimates \eqref{stime m3 alpha} and by Lemma \ref{stac}, one obtains that for 
$0 < \zeta \ll \rho$
\begin{equation}\label{stime lambda 3 + beta}
|\lambda_3^+ - \lambda_3|^{\Omega} \lesssim \delta, \quad | \beta |_{\rho - \zeta}^{\Omega} \lesssim  
 e^{\tC_0 \zeta^{- \mu}} \delta\,.
\end{equation}
By Lemma \ref{lemma diffeo inverso} and \eqref{ruota} we see that
$\f \mapsto \f + \omega \beta(\f)$ is invertible and the inverse is given by 
$\vartheta \mapsto \vartheta + \omega \widetilde \beta(\vartheta)$ with 
\begin{equation}\label{stima alpha tilde alpha sn}
\widetilde \beta \in {\cal H}(\T^\infty_{\rho - 2\zeta}), \quad | \widetilde \beta |_{\rho - 2\zeta }^{\Omega} \lesssim 
e^{\tC_0 \zeta^{- \mu}} \delta\,. 
\end{equation}

The inverse of the operator ${\cal T}_2$ is then given by
\begin{equation}\label{def cal B inverso}
{\cal T}_2^{- 1} u(\vartheta, x) = u(\vartheta + \omega \widetilde \beta(\vartheta), x)\,. 
\end{equation}
so that
\begin{equation}\label{coniugio B cal L (1)}
\begin{aligned}
{\cal T}_2^{- 1} {\cal L}^{(1)} {\cal T}_2 &= {\cal T}_2^{- 1}
\big( 1 + \omega \cdot \partial_\f \beta \big) \omega \cdot \partial_\vartheta + 
{\cal T}_2^{- 1}(m_3) \partial_x^3 + {\cal T}_2^{- 1}(b_1) \partial_x + {\cal T}_2^{- 1}(b_0)\\
&=:\frac{1}{\mathtt r} {\cal L}^{(2)}
\end{aligned}
\end{equation}
where 
\begin{equation}\label{cal L (2) finale}
\begin{aligned}
& {\cal L}^{(2)} := \omega \cdot \partial_\vartheta + \lambda_3^+ \partial_x^3 + 
c_1(\vartheta, x) \partial_x + c_0(\vartheta, x)\,, \\
& {\mathtt r} := \frac{1}{{\cal T}_2^{- 1}\big(1 + \omega \cdot \partial_\f \be \big)}
 \stackrel{\eqref{eq omologica riparametrizzazione tempo}}{=}
  \frac{\lambda_3^+}{{\cal T}_2^{- 1}(m_3)}\,, \\
& c_i  := {\mathtt r}{{\cal T}_2^{- 1}(b_i)}, \quad i = 1, 0\,.  
\end{aligned}
\end{equation}

Therefore by the estimates \eqref{stime m3 alpha}, \eqref{stime lambda 3 + beta}, 
\eqref{stima alpha tilde alpha sn} and by applying Corollary 
\ref{usiamoquesto}, Lemma \ref{moser type lemma},
and \eqref{ruota}, one gets  
\begin{equation}\label{voila}
\begin{aligned}
&| {\mathtt r} - 1 |_{\rho - \zeta}^{\Omega} \lesssim \delta\,\\
&| c_i - a_i |_{\rho - \zeta}^{\Omega} 
\,{\lesssim}\, \zeta^{- \tau}e^{\tC_0 \zeta^{- \mu }} \delta , \quad i = 0,1\,.
\end{aligned}
\end{equation}

%%%%%%%%%%%%%%%%%%%%%%%%%%%%%%%%%%%%%%%%%%%%%%%%%
\subsection{Time dependent traslation of the space variable}
%%%%%%%%%%%%%%%%%%%%%%%%%%%%%%%%%%%%%%%%%%%%%%%%%

Let $p : \T^\infty_{\rho-2\ze} \to \R$ be an analytic function to be determined and let 
\begin{equation}\label{def cal T3}
{\cal T}_3 u(\f, x) := u(\f, x + p(\f)), \quad \text{with inverse} \quad {\cal T}_3^{- 1} u(\f, y) = u(\f, y - p(\f))\,. 
\end{equation}

Computing explicitly
\begin{equation}\label{cal L (3)}
\begin{aligned}
{\cal L}^{(3)} &:= {\cal T}_3^{- 1} {\cal L}^{(2)} {\cal T}_3 =  
 \omega \cdot \partial_\f + \lambda_3^+ \partial_x^3 + a_1^+(\f, x) \partial_x + a_0^+(\f, x)\,, \\
& a_1^+  := \omega \cdot \partial_\f p + {\cal T}_3^{- 1}(c_1)\,, \quad a_0^+ := {\cal T}_3^{- 1}(c_0)\,,
\end{aligned}
\end{equation}
and by \eqref{prop mediaMic} one has  
\begin{equation}\label{media p rn}
\begin{aligned}
\frac{1}{2 \pi} \int_\T {\cal T}_3^{- 1}(c_1)(\f, y)\, d y & =
 \frac{1}{2 \pi} \int_\T c_1(\f, x)\, d x \\
 &= \frac{1}{2 \pi} \int_\T a_1(\f, x)\, d x +
  \frac{1}{2 \pi} \int_\T (c_1 - a_1)(\f, x)\, d x   \\
& = \lambda_1 + \frac{1}{2 \pi} \int_\T (c_1 - a_1)(\f, x)\, d x\,.
\end{aligned}
\end{equation}

We want to choose $p(\f)$ in such a way that the $x$-average of $d_1$ is constant.
To this purpose we define 
\begin{equation}\label{def p}
p(\f) := (\omega \cdot \partial_\f)^{- 1}\Big[ \langle (c_1 - a_1)\rangle_{\f, x} -
 \frac{1}{2 \pi} \int_\T (c_1 - a_1)(\f, x)\, d x \Big]
\end{equation}
where for any $a : \T^\infty_\sigma \times \T_\sigma \to \C$, $\langle a \rangle_{\f, x}$ is defined by 
$$
\langle a \rangle_{\f, x} := \frac{1}{(2 \pi)} \int_\T \int_{\T^\infty} a(\f, x)\, d \f\, d x
$$
(recall the definition \ref{media toro infinito dimensionale}). By  \eqref{voila} and  Lemma \ref{stac} one gets
\begin{equation}\label{stima p}
| p |_{\rho - 2\zeta}^\Omega \lesssim 
\zeta^{- \tau} e^{2 \tC_0 \zeta^{- \mu }} \delta \stackrel{\eqref{ruota}}{\ll}\ze \,.
\end{equation}

Moreover  
\begin{equation}\label{cal L (3) finale}
 \lambda_1^+ := \frac{1}{2 \pi} \int_\T d_1(\f, x)\, d x = \lambda_1 + \langle (c_1-a_1) \rangle_{\f, x}\,. 
\end{equation}

Finally using \eqref{voila}, \eqref{cambio variabile meno identita} 
(with $\Phi_\alpha = {\cal T}_3^{- 1}$), \eqref{stima p}, one gets 
\begin{equation}\label{di - ai}
| a^+_i - a_i |_{\rho - 2\zeta}^\Omega \lesssim \zeta^{- \tau'}e^{2\tC_0 \zeta^{- \mu }} \delta \,,
\end{equation}
for some $\tau'>0$.

%%%%%%%%%%%%%%%%%%%%%%%%%%%%%%%%%%%%%%%%%%%%%%%%%
\subsection{Conclusion of the proof}
%%%%%%%%%%%%%%%%%%%%%%%%%%%%%%%%%%%%%%%%%%%%%%%%%

We start by noting that $\TT:=\TT_3\circ \TT_2\circ\TT_1$ has the form \eqref{nonsicak}
with $p^{(n+1)}=p$, $\be^{(n+1)}=\be$ and $\x^{(n+1)}(\f,x) = \al(\f+\om\be(\f),x+p(\f))$.
Hence, setting $ {\mathtt r}:={\mathtt r}_{n+1}$, $\rho := s_n-\s_n$, $\de := \s_n^{-4}\e_n$,
 $\de_0:=2\e_0$ and $\ze:= \s_n$ we denote
 \[
 1+A_{n+1} = \la_3^+,\quad, B_{n+1}(\f,x):= a_1^{+}(\f,x),\quad
 C_{n+1}=a_0^+(\f,x),
 \] 
 and thus Proposition \ref{comecaz} follows.
 \EP

%%%%%%%%%%%%%%%%%%%%%%%%%%%%%%%%%%%%%%%%%%%%%%%%%
%%%%%%%%%%%%%%%%%%%%%%%%%%%%%%%%%%%%%%%%%%%%%%%%%
\section{Proof of Proposition \ref{ciaraggio}}\label{eppoi}
\zerarcounters
%%%%%%%%%%%%%%%%%%%%%%%%%%%%%%%%%%%%%%%%%%%%%%%%%
%%%%%%%%%%%%%%%%%%%%%%%%%%%%%%%%%%%%%%%%%%%%%%%%%

In order to prove Proposition \ref{ciaraggio},  we start by considering
 a linear Hamiltonian operator
 defined for $\omega \in {\cal O} \subseteq \mathtt D_\gamma$
of the form 
\begin{equation}\label{espansione di finale}
{\cal L} = {\cal L}(\lambda_3, a_1, a_0) := \omega \cdot \partial_\f + \lambda_3 \partial_x^3 + 
a_1(\f, x) \partial_x + a_0(\f, x)\,.\\
\end{equation}

We want to show that, for any choice of the coefficients $\la_3,a_1,a_0$ satisfying some hypotheses
(see below), it is possible to reduce $\LL$ to constant coefficients.
Moreover we want to show
that such reduction is ``Lipshitz'' w.r.t. the parameters $\la_3,a_1,a_0$, in a sense that will be
clarified below.

Regarding the coefficients, we need to require that
\begin{equation}\label{risistemo}
\begin{aligned}
& a_i :=  \sum_{k = 0}^{m} a_i^{(k)}, \quad | a_i^{(k)} |_{\rho_k}^{\cal O} \lesssim 
\delta_k, \quad \forall k = 0, \ldots, m\,,\ \ i=0,1,  \\
& |\lambda_3 -1|^{\cal O}\lesssim \delta_0 \,, \\
& \lambda_1 \equiv \lambda_1(a_1) = \sum_{k = 0}^{m} \lambda_1^{(k)},
 \quad \lambda_1^{(k)} := \frac{1}{2 \pi} \int_\T a_1^{(k)}(\f, x)\, d x = {\rm const}\,. 
\end{aligned}
\end{equation}
for some
$0<\ldots <\rho_m <  \ldots < \rho_0$ and $0<\ldots\ll\delta_m  \ll \ldots \ll \delta_0\ll1$
so that there is a third sequence $\zeta_i$ such that $0 < \zeta_i < \rho_i$ and 
\begin{equation}\label{ipotesi sommabilita ei}
\sum_{i \ge0} \zeta_i^{- \tau} { e}^{C \zeta_i^{- \mu}} \delta_i \lesssim \de_0 \,,
\end{equation}
for some $\tau,C>0$.

%%%%%%%%%%%%%%%%%%%%%%%%%%%%%%%%%%%%%%%%%%%%%%%%%
\subsection{Reduction of the first order term}
%%%%%%%%%%%%%%%%%%%%%%%%%%%%%%%%%%%%%%%%%%%%%%%%%

We consider an operator $\LL$ of the form \eqref{espansione di finale} 
satisfying the hypotheses above. We start by showing that it is possible to
reduce it
to constant coefficients up to a bounded reminder, and that such reduction is 
``Lipshitz'' w.r.t. the parameters $\la_3,a_1,a_0$.

\begin{lemma}\label{riduzione ordine 1}
There exists a symplectic invertible operator ${\cal M} = {\rm exp}({\cal G})$, 
with ${\cal G} \equiv {\cal G}(\lambda_3, a_1)$ and an 
operator ${\cal R}_0 \equiv {\cal R}_0(\lambda_3, a_1, a_0)$ satisfying 
\begin{equation}\label{prop cal G cal R0}
\begin{aligned}
& {\cal G } = \sum_{i = 0}^{m } {\cal G}^{(i)}\,, \quad \|{\cal G}^{(i)} \|_{\rho_i, - 1}^{\cal O} \lesssim \delta_i\,, \\
& {\cal R}_0 = \sum_{i = 0}^{m } {\cal R}_0^{(i)}, \quad \| {\cal R}_0^{(i)} \|_{\rho_i - \zeta_i}^{\cal O} 
\lesssim \zeta_i^{- \tau} e^{C \zeta_i^{- \mu}} \delta_i
\end{aligned}
\end{equation}
for some $C , \tau \gg 1$, such that 
\begin{equation}\label{coniugio cal L cal L 0}
{\cal L}_0 := {\cal M}^{- 1} {\cal L} {\cal M} = \omega \cdot \partial_\f + \lambda_3 \partial_x^3 +
 \lambda_1 \partial_x + {\cal R}_0\,. 
\end{equation}
\end{lemma}

\prova
We look for $\calG$ of the form 
\[
\calG = \pi_0^\bot g(\f, x) \partial_x^{- 1}
\]
and we choose the function $g (\f, x)$ where $g = g(\lambda_3,   a_1)$ in order to solve 
\begin{equation}\label{eq omologica riduzione partial x}
\begin{aligned}
& 3 \lambda_3 \partial_x g(\f, x) + a_1(\f, x) = \lambda_1\,.
\end{aligned}
\end{equation}
By \eqref{risistemo}, one obtains that 
\begin{equation}\label{def g a}
g := \frac{1}{3 \lambda_3} \partial_x^{- 1}\Big[ \lambda_1 - a_1\Big]
\end{equation}
and therefore
\begin{equation}\label{decomposizione g}
\begin{aligned}
& g = \sum_{i = 0}^{m} g_i\,, \quad  g_i := \frac{1}{3 \lambda_3} \partial_x^{- 1}
\Big[ \lambda_1^{(i)} - a_1^{(i)}\Big]\,, \\
& | g_i |_{\rho_i}^\calO \lesssim  \delta_i, \quad i = 0, \ldots, m\,. \\
\end{aligned}
\end{equation}

 Of course we can also write the operator ${\cal G} := \pi_0^\bot g(\f, x) \partial_x^{- 1} 
 = \sum_{i = 0}^{m } {\cal G}_i$ where ${\cal G}_i := \pi_0^\bot g_i(\f, x) \partial_x^{- 1}$ and  one has 
\begin{equation}\label{stima cal Gi}
\begin{aligned}
& \|{\cal G}_i \|_{\rho_i, - 1}^\calO \lesssim  \delta_i, \quad i = 0, \ldots,  m \,.
\end{aligned}
\end{equation}

Again by \eqref{risistemo}, defining ${\cal P} := a_1 \partial_x + a_0$, one has that 
${\cal P} = \sum_{i = 0}^{m} {\cal P}_i$, where ${\cal P}_i := a_1^{(i)} \partial_x + a_0^{(i)}$ satisfies   
\begin{equation}\label{stime cal Pi ordine 1}
\|{\cal P}_i \|_{\rho_i , 1}^\calO \lesssim  \delta_i.
\end{equation}

Therefore
\begin{equation}\label{operator cal L (4)}
\begin{aligned}
{\cal L}_0 &=  {\cal M}^{- 1} {\cal L} {\cal M} = e^{- {\cal G}}  \omega \cdot \partial_\f   e^{\cal G} +
 \lambda_3 e^{- {\cal G}} \partial_x^3 e^{\cal G} + e^{- {\cal G}}  {\cal P}e^{\cal G}  \\
& = \omega \cdot \partial_\f + \lambda_3 \partial_x^3 + \Big(3 \lambda_3 g_x + a_1\Big) \partial_x + {\cal R}_0 \\
& \stackrel{\eqref{eq omologica riduzione partial x}}{=} 
\omega \cdot \partial_\f + \lambda_3 \partial_x^3 + \lambda_1 \partial_x + {\cal R}_0
\end{aligned}
\end{equation} 
where %the remainder ${\cal R}_0 \equiv {\cal R}_0(\lambda_3, a_1, a_0)$ is defined by 
\begin{equation}\label{definizione resto step dx}
\begin{aligned}
{\cal R}_0 & := \Big( e^{- {\cal G}}  \omega \cdot \partial_\f   e^{\cal G} - 
\omega \cdot \partial_\f \Big) + \lambda_3 \Big( e^{- {\cal G}} \partial_x^3 e^{\cal G} -
 \partial_x^3 - 3 g_x \partial_x  \Big)  + \Big( e^{-{\cal G}} {\cal P} e^{\cal G} - {\cal P} \Big) + a_0\,. \\
\end{aligned}
\end{equation}

Then \eqref{ipotesi sommabilita ei}, \eqref{stima cal Gi}, \eqref{stime cal Pi ordine 1} 
guarantee that the hypotheses of Lemmata 
\ref{lemma commutatore R ord 1}-\ref{commutatore epsilon k sigma k} are verified. 
Hence, we apply Lemma \ref{lemma commutatore R ord 1}-$(ii)$ to expand the operator 
$e^{-{\cal G}} {\cal P} e^{\cal G} - {\cal P}$, Lemma \ref{commutatore epsilon k sigma k}-$(ii)$ 
to expand $e^{- {\cal G}} \partial_x^3 e^{\cal G} - \partial_x^3 - 3 g_x \partial_x$ and 
Lemma \ref{commutatore epsilon k sigma k}-$(iii)$ to expand 
$e^{- {\cal G}}  \omega \cdot \partial_\f   e^{\cal G} - \omega \cdot \partial_\f$. 
The expansion of the multiplication operator $a_0$ is already provided by \eqref{risistemo}. 
Hence, one obtains that there exist $C, \tau \gg 1$  such that \eqref{prop cal G cal R0}
 is satisfied.
\EP

We now consider a ``small modification'' of the operator $\LL$ in the following sense. 
We consider an operator
\begin{equation}\label{coniugio cal L0 +}
{\cal L}^+=\LL(\la_3^+,a_1^+,a_0^+) : = \omega \cdot \partial_\f + \lambda_3^+ \partial_x^3 + a_1^+(\f,x)
\partial_x +a_0^+(\f,x)
\end{equation}
with
\begin{equation}\label{ipotesi variazioni lemma partial x}
 \frac{1}{2\pi} \int_\T a_1^+(\f, x)\, d x =:\la_1^+= {\rm const},
  \quad |a_i^+ - a_i |_{\rho_{m+1}}\,,\, |\lambda_3^+ - \lambda_3|\lesssim \delta_{m+1}\,. 
\end{equation}

Of course we can apply Lemma \ref{riduzione ordine 1} and conjugate $\LL^+$
to 
\begin{equation}\label{coniugio cal L0 ++}
{\cal L}_0^+ : = \omega \cdot \partial_\f + \lambda_3^+ \partial_x^3 + \lambda_1^+ \partial_x + 
{\cal R}_0^+
\end{equation}
with $\RR_0^+$ a bounded operator. We want to show that $\LL_0^+$ is ``close'' to $\LL_0$,
namely the following result.

\begin{lemma}\label{lemma variazione parte partial x}
One has
\begin{equation}\label{stima G - G + R0 - R0 +}
\begin{aligned}
|\lambda_1^+ - &\lambda_1| \lesssim \delta_{m+1}\,,\qquad
%\quad  \| {\cal G}^+ - {\cal G} \|_{\rho_{m+1}, - 1} 
%\lesssim \delta_{m+1}, \\
%&
 \| {\cal R}_0^+ - {\cal R}_0 \|_{\rho_{m+1} - \zeta_{m+1}} \lesssim
 \zeta_{m+1}^{- \tau}e^{C \zeta_{m+1}^{- \mu }} \delta_{m+1}\,. 
 \end{aligned}
\end{equation}
\end{lemma}

\prova
The first bound follows trivially from \eqref{ipotesi variazioni lemma partial x}.
Regarding the second bound one can reason as follows.
As in Lemma \ref{riduzione ordine 1}, er can define 
${\cal G}^+ := \pi_0^\bot g^+(\f, x) \partial_x^{- 1}$ with 
\begin{equation}\label{def g +}
 g^+ := \frac{1}{3 \lambda_3^+} \partial_x^{- 1}\Big[ \lambda_1^+ - a_1^+\Big]
\end{equation}
so that
\begin{equation}\label{decomposizione g +}
\begin{aligned}
 \| {\cal G}^+ - {\cal G} \|_{\rho_{m+1}, -1} \lesssim \delta_{m+1}\,. 
\end{aligned}
\end{equation}

Defining ${\cal P}^+ := a_1^+ \partial_x + a_0^+$ and recalling that 
${\cal P} := a_1 \partial_x + a_0$, by \eqref{ipotesi variazioni lemma partial x}, one gets 
\begin{equation}\label{stime cal Pi ordine 1 +}
 \|{\cal P}^+ - {\cal P} \|_{\rho_{m+1}, 1} \lesssim \delta_{m+1}\,. 
\end{equation}

The estimate on ${\cal R}_0^+ - {\cal R}_0$ follows by applying Lemmata \ref{variazione 1}, \ref{variazione 2},
 and by the estimates \eqref{ipotesi variazioni lemma partial x}, \eqref{stime cal Pi ordine 1 +}, 
 \eqref{decomposizione g +}.  
\EP

%%%%%%%%%%%%%%%%%%%%%%%%%%%%%%%%%%%%%%%%%%%%%%%%%
\subsection{Reducibility}\label{pelle}
%%%%%%%%%%%%%%%%%%%%%%%%%%%%%%%%%%%%%%%%%%%%%%%%%

We now consider an operator $\LL_0$ of the form
\begin{equation}\label{cia}
{\cal L}_0 \equiv {\cal L}_0(\lambda_1, \lambda_3, {\cal P}_0) 
:= \omega \cdot \partial_\f + {\cal D}_0 + {\cal P}_0
\end{equation}
with $\calP_0$ a bounded operator and
\begin{equation}\label{cal D (0)}
 {\cal D}_0 \equiv{\cal D}_0(\lambda_1, \lambda_3) := 
 \ii\, {\rm diag}_{j \in \Z \setminus \{  0 \}} \Omega_{0}(j)\,, \qquad
 \Omega_{ 0}(j) := - \lambda_3 j^3 +  \lambda_1 j, \quad j \in \Z \setminus \{ 0 \}\,, 
\end{equation}
and we show that, under some smallness conditions specified below
it is possible to reduce it to constant coefficients, and that the reduction
is ``Lipschitz'' w.r.t. the parameters $\lambda_1, \lambda_3, {\cal P}_0$.

In order to do so, we introduce three sequences $ 0<\ldots< \rho_{m } < \ldots < \rho_0$, 
$0<\ldots \ll \delta_{m } \ll \ldots \ll \delta_0$ and $1\ll N_0 \ll N_1 \ll \cdots$
 and we assume that setting
{$\Delta_i=\rho_i-\rho_{i+1}$} one has
\begin{equation}\label{lostrozzo}
\sum_{i \ge0} \Delta_i^{- \tau} { e}^{C \Delta_i^{- \mu}} \delta_i \lesssim\de_0\,,
\end{equation} 
\begin{equation}\label{dafare}
	 e^{- N_k\Delta_k}  \delta_k+
	 e^{C \Delta_k^{- \mu}} \delta_k^2 
	 \ll  2^{-k}  \delta_{k + 1}\,,
\end{equation}
 \begin{equation}\label{dafare2}
 \delta_k \ll (1+N_k)^{-C N_k^{\frac{1}{1+\h}}}	\end{equation}
 and 
\begin{equation}\label{gufo}
\begin{aligned}
& |\lambda_3-1|^{\cal O}, |\lambda_1|^{\cal O} \le \delta_0, \\
& {\cal P}_0 :=  \sum_{i = 0}^{m} {\cal P}_0^{(i)}\,, \quad  \| {\cal P}_0^{(i)} \|_{\rho_i }^{\cal O} 
\le   \delta_i ,\quad i = 0, \ldots, m\,,
\end{aligned}
\end{equation}
for some $\tau,C>0$.

We have the following result.

\begin{lemma}\label{lemma iterativo riducibilita}
Fix $\g\in[\g_0/2,2\g_0]$.
For $k=0,\ldots,m$ there is a sequence of sets $\calE_k\subseteq\calE_{k-1}$ and a sequence of
symplectic maps $\Phi_k$ defined for $\om\in\calE_{k+1}$ such that setting 
$\LL_0$ as in \eqref{cia} and for $k\ge1$,
\begin{equation}\label{coniugio riducibilita mortacci}
\LL_{k}:= \Phi_{k - 1}^{- 1} {\cal L}_{k - 1} \Phi_{k - 1},
\end{equation}
one has the following.

\begin{enumerate}

\item 
$\LL_k$ is of the form
\begin{equation}\label{cal L (k) red porchi}
{\cal L}_k := \omega \cdot \partial_\f + {\cal D}_k + {\cal P}_k
\end{equation}

where

\begin{itemize}

\item 
The operator $\DD_k$ is of the form
	\begin{equation}\label{cal D (k) rid}
	\begin{aligned}
	& {\cal D}_k := {\rm diag}_{j \in \Z \setminus \{ 0 \}} \Omega_k(j), \quad \Omega_k(j) = \Omega_0 (j) + r_k(j)
	\end{aligned}
	\end{equation}
	with $r_0(j)=0$ 
	and for $k\ge1$, $r_k(j)$ is defined for $\omega \in {\calE}_0=\calO$ and satisfies
	\begin{equation}\label{Omega j (k + 1) (k)}
	\begin{aligned}
	&  \sup_{j \in \Z \setminus \{ 0 \}}|r_k(j) - r_{k - 1}(j)|^{\cal O} \le  \delta_{k - 1} \sum_{i=1}^{k - 1}2^{-i} \,. 
	 \end{aligned}
	\end{equation}
	
\item
The operator ${\cal P}_k$ is such that
	\begin{equation}\label{operator cal P (k)}
	\begin{aligned}
	& \text{for} \quad 0 \leq k \leq m, \quad 
	{\cal P}_k = \sum_{i = k}^{m} {\cal P}_k^{(i)}\,,\quad 
	 \|{\cal P}_k^{(i)}\|_{\rho_i}^{{\calE}_k} \le  \delta_i \sum_{j=1}^{k} 2^{-j}, \quad \forall i = k, \ldots, m\,. \\
	%& \|{\cal P}_n \|_{\rho_n}^{{\calE}_n} \leq C_*^n  \delta_n \\
	\end{aligned}
	\end{equation}

\end{itemize}

\item
One has $\Phi_{k - 1} := {\rm exp}(\Psi_{k - 1})$, such that
	\begin{equation}\label{stima Psi k - 1 lemma}
	\begin{aligned}
	 \|\Psi_{k - 1} \|_{\rho_k}^{{\calE}_k} \lesssim 
	e^{C \Delta_{k-1}^{- \mu}}
	\|{\cal P}_{k - 1}^{(k - 1)}\|_{\rho_{k - 1}}^{{\calE}_{k - 1}} 
	 \lesssim 
	e^{C \Delta_{k-1}^{- \mu}} \delta_{k - 1} \\
	\end{aligned}
	\end{equation}
	
\item The sets $\calE_k$ are defined as
	\begin{equation}\label{def cal O (k)}
	\begin{aligned}
	{\calE}_k  := \Big\{\omega \in {\calE}_{k - 1} \ : \
	&|\omega \cdot \ell + \Omega_{k - 1}(j) - \Omega_{ k - 1}(j')|  \geq \frac{\gamma |j^3 - j'^3|}{\mathtt d(\ell)},  \\
	& \quad \forall (\ell , j, j') \neq (0, j, j), \quad |\ell|_\eta \leq N_{k - 1} \Big\}.
		\end{aligned}
	\end{equation}

\end{enumerate}
\end{lemma}

\prova
The statement is trivial for $k=0$ so we assume
it to hold up to $k<m$ and let us prove it for $k + 1$.
For any $\Phi_k := {\rm exp}(\Psi_k)$ one has 
\begin{equation}\label{def cal L (k + 1)}
\begin{aligned}
& {\cal L}_{k + 1}= \Phi_k^{- 1} {\cal L}_k \Phi_k
 = \omega \cdot \partial_\f + {\cal D}_k + \omega \cdot \partial_\f \Psi_k + 
[{\cal D}_k, \Psi_k]+ \Pi_{N_k}{\cal P}^{(k)}_k + {\cal P}_{k + 1} 
\end{aligned}
\end{equation}
where the operator ${\cal P}_{k + 1}$ is defined by 
\begin{equation}\label{def cal P (k + 1)}
\begin{aligned} 
{\cal P}_{k + 1} & := \Pi_{N_k}^\bot {\cal P}_k^{(k)} + 
\sum_{p \geq 2} \frac{{\rm Ad}_{\Psi_k}^p(\omega \cdot \partial_\f + {\cal D}_k)}{p!} +
 \sum_{i = k + 1}^{m} e^{- \Psi_k} {\cal P}^{(i)}_k e^{\Psi_k} + 
 \sum_{p \geq 1} \frac{{\rm Ad}^p_{\Psi_k}({\cal P}^{(k)}_k)}{p!}\,.
\end{aligned}
\end{equation}

Then we choose $\Psi_k$ in such a way that 
\begin{equation}\label{eq omologica Psi k}
\begin{aligned}
& \omega \cdot \partial_\f \Psi_k + [{\cal D}_k, \Psi_k]+ \Pi_{N_k} {\cal P}^{(k)}_k = {\cal Z}_k\,, \\
& {\cal Z}_k := {\rm diag}_{j \in \Z \setminus \{ 0 \}}  ({\cal P}^{(k)}_k)_j^j(0)\,,
\end{aligned}
\end{equation}
namely for $\omega \in {\calE}_{k + 1}$ we set
\begin{equation}\label{def matrice Psi k}
( \Psi_k)_j^{j'}(\ell) :=  \begin{cases}
\dfrac{({\cal P}_k^{(k)})_j^{j'}(\ell)}{\ii \big(\omega \cdot \ell + \Omega_{k}(j) - \Omega_k(j') \big)}, 
 \quad \forall (\ell, j, j') \neq (0, j, j), \quad |\ell|_\eta \leq N_k\,, \\
0 \quad \quad \quad  \text{otherwise.}
	\end{cases}
	\end{equation}
	
Therefore, 
\begin{equation}\label{stima coefficienti matrice Psi k}
 | ( \Psi_k)_j^{j'}(\ell)| \lesssim  \mathtt d(\ell) |({\cal P}_k^{(k)})_j^{j'}(\ell)|\,, \quad \forall \omega \in {\calE}_{k+1}\,. 
\end{equation}
and by applying Lemma \ref{bound per stima di Cauchy}, using the induction estimate 
\eqref{operator cal P (k)}, one obtains
 \begin{equation}\label{stima Psi (k)}
 \|\Psi_k \|_{\rho_k - \zeta}^{{\calE}_{k+1} }
 \lesssim  e^{C \zeta^{- \mu}} \|{\cal P}_k^{(k)} \|_{\rho_k}^{{\calE}_{k} }
 \stackrel{\eqref{operator cal P (k)}}{\lesssim} e^{C \zeta^{- \mu}}  \delta_k\,, \\
\end{equation}
for any $\zeta<\rho_k$.

We now define the diagonal part ${\cal D}_{k + 1}$.

For any $j \in \Z \setminus \{ 0 \}$ and any 
$\omega \in {\calE}_k$  one has
 $|({\cal P}_k^{(k)})_j^j(0)| \lesssim \|{\cal P}_k^{(k)}\|_{\rho_k}^{{\calE}_{k} }
 \stackrel{\eqref{operator cal P (k)}}{\leq} \delta_k $. 
 The Hamiltonian structure guarantees that ${\cal P}_k^{(k)}(0)_j^j$ is purely imaginary 
 and by the Kiszbraun Theorem there exists a Lipschitz extension 
 $\omega \in {\cal O} \to \ii z_k(j)$ (with $z_k(j)$ real) 
 of this function satisfying the bound $|z_k(j)|^\calO \lesssim  \delta_k$. Then, we define
	\begin{equation}\label{def nuova diagonale (k + 1)}
	\begin{aligned}
	& {\cal D}_{k + 1} :=  {\rm diag}_{j \in \Z \setminus \{ 0 \}} \Omega_{k + 1}(j)\,, \\
	& \Omega_{k + 1}(j) := \Omega_k(j) + z_k(j) = 
	\Omega_0(j) + r_{k + 1}(j), \quad \forall j \in \Z \setminus \{ 0 \}\,, \\
	& r_{k + 1}(j) := r_k(j) + z_k(j) \\
	\end{aligned}
	\end{equation}
	and one has 
	\begin{equation}\label{stima Omega j (k + 1) (k)}
	|r_{k + 1}(j) - r_k(j)|^{\calO}= |z_k(j)|^{\calO} \leq \|{\cal P}_k^{(k)} \|_{\rho_k}^{{\calE}_{k} } 
	\stackrel{\eqref{operator cal P (k)}}{\leq}  \delta_k \sum_{j=1}^k 2^{-j} 
	\end{equation}
	which is the estimate \eqref{Omega j (k + 1) (k)} at the step $k + 1$.

We now estimate the remainder ${\cal P}_{k + 1}$ in \eqref{def cal P (k + 1)}. Using
\eqref{eq omologica Psi k} we see that 
\begin{equation}\label{def cal P (k + 1) b}
\begin{aligned} 
{\cal P}_{k + 1}& = \Pi_{N_k}^\bot {\cal P}_k^{(k)} + 
\sum_{p \geq 2} \frac{{\rm Ad}_{\Psi_k}^{p - 1}({\cal Z}_k - \Pi_{N_k}{\cal P}_k^{(k)})}{p!} +
 \sum_{i = k + 1}^{m} e^{- \Psi_k} {\cal P}^{(i)}_k e^{\Psi_k} + 
 \sum_{p \geq 1} \frac{{\rm Ad}^p_{\Psi_k}({\cal P}^{(k)}_k)}{p!}\,.  \\
\end{aligned}
\end{equation}

Denote 
	\begin{equation}\label{decomposizione cal P (k + 1)}
	\begin{aligned}
	& {\cal P}_{k + 1}= \sum_{i = k + 1}^{m} {\cal P}_{k + 1}^{(i)} \quad \text{where} \\
	& {\cal P}_{k + 1}^{(k + 1)} := \Pi_{N_k}^\bot {\cal P}_k^{(k)} + 
	 \sum_{p \geq 2} \frac{{\rm Ad}_{\Psi_k}^{p - 1}({\cal Z}_k - \Pi_{N_k}{\cal P}_k^{(k)})}{p !} 
	 + e^{- \Psi_k} {\cal P}^{(k + 1)}_{k } e^{\Psi_k} + 
	 \sum_{p \geq 1} \frac{{\rm Ad}^p_{\Psi_k}({\cal P}^{(k)}_k)}{p!}\,, \\
	& {\cal P}_{k + 1}^{(i)} := e^{- \Psi_k} {\cal P}^{(i)}_k e^{\Psi_k}, \quad i = k + 2, \ldots, m\,. 
	\end{aligned}
	\end{equation}

\noindent
	{\sc Estimate of ${\cal P}_{k + 1}^{(i)}$, $i = k + 2, \ldots, m$.} By the induction estimate, one has
	\begin{equation}\label{stima cal P (k + 1) i geq 2}
	\begin{aligned}
	\|e^{- \Psi_k}{\cal P}_k^{(i)} e^{\Psi_k} \|_{\rho_i}^{{\calE}_{k+1} } &\le
	\|{\cal P}_k^{(i)} \|_{\rho_i}^{{\calE}_{k} } +	
	\|{\cal P}_k^{(i)}-e^{- \Psi_k}{\cal P}_k^{(i)} e^{\Psi_k} \|_{\rho_i}^{{\calE}_{k+1} }  \\
	& \lesssim   \delta_i \sum_{j=1}^{k}2^{-j} +  \| \Psi_k\|_{\rho_i}^{{\calE}_{k+1} }
	\| \calP^{(i)}_k\|_{\rho_i}^{{\calE}_{k+1} }  \stackrel{\eqref{dafare}}{\lesssim} \delta_i \sum_{j=1}^{k+1}2^{-j} \,.
	 \end{aligned}
	\end{equation}

\noindent	
{\sc Estimate of ${\cal P}_{k + 1}^{(k + 1)}$.} We estimate separately the four terms in the definition of 
${\cal P}_{k + 1}^{(k + 1)}$  in \eqref{decomposizione cal P (k + 1)}.
 By Lemma \ref{proprieta norma sigma riducibilita vphi x}-$(ii)$, one has  
	\begin{equation}\label{stima Pi Nk bot P kk}
	\|\Pi_{N_k}^\bot {\cal P}_k^{(k)} \|_{\rho_{k + 1}}^{\calE_k} \lesssim 
	e^{- N_k\Delta_k} \|{\cal P}_k^{(k)} \|_{\rho_k}^{\calE_k} \lesssim 
	e^{- N_k\Delta_k}  \delta_k\,. 
	\end{equation}
	
By applying \eqref{stima banale Ad n} and the estimate of Lemma 
\ref{proprieta norma sigma riducibilita vphi x}-$(iii)$,  one obtains 
	\begin{equation}\label{primo pezzo cal P 1 (1)}
	\begin{aligned}
	 \Big\| \sum_{p \geq 2} 
	&\frac{{\rm Ad}_{\Psi_k}^{p - 1}({\cal Z}_k - \Pi_{N_k} {\cal P}_k^{(k)})}{p!}  \Big\|_{\rho_{k + 1}}^{\calE_{k+1}} 
	\leq \sum_{p \geq 2} \frac{C^{p - 1}}{p!} (\|\Psi_k \|_{\rho_{k + 1}}^{\calE_{k+1}})^{p - 1} 
	\|{\cal P}_k^{(k)} \|_{\rho_k}^{\calE_{k}}  \\
	& {\lesssim}
	 \|\Psi_k \|_{\rho_{k + 1}}^{\calE_{k+1}}  \|{\cal P}_k^{(k)} \|_{\rho_k}^{\calE_{k}}  \lesssim  
	 e^{C \Delta_k^{- \mu}}   \delta_k^2 
	\end{aligned}
	\end{equation}
	and similarly 
	\begin{equation}\label{stima terzo pezzo cal P k + 1 (k + 1)}
	\Big\| \sum_{m \geq 1} \frac{{\rm Ad}^m_{\Psi_k}({\cal P}^{(k)}_k)}{m!} \Big\|_{\rho_{k + 1}}^{\calE_{k+1}}
	 \lesssim e^{C \Delta_k^{- \mu}}  \delta_k^2\,. 
	\end{equation}

In conclusion we obtained	
	\begin{equation}\label{stima finale cal P k + 1 (k + 1)}
	\|{\cal P}_{k + 1}^{(k + 1)}  \|_{\rho_{k + 1}}^{\calE_{k+1}} 
	\le C' e^{- N_k\Delta_k}   \delta_k +
  C'e^{C\Delta_k^{- \mu }}  \delta_k^2 +  \delta_{k + 1}\sum_{j=1}^{k}2^{-j}
	\end{equation}
where $C'$ is an appropriate constant and the last summand is a bound for the term 
$e^{- \Psi_k} {\cal P}^{(k + 1)}_{k } e^{\Psi_k}$, which
can be obtained reasoning as in \eqref{stima cal P (k + 1) i geq 2}.
Thus we obtain  
\begin{equation}\label{thus}
 \|{\cal P}_{k+1}^{(k+1)}\|_{\rho_{k+1}}^{{\calE}_{k+1}} \le  \delta_{k + 1} \sum_{j=1}^{k+1} 2^{-j}
\end{equation}
provided 
	$$
	 C' e^{- N_k\Delta_k} \delta_k+ 
	C'  e^{C\Delta_k^{- \mu}}
	  \delta_k^2 +
	   \delta_{k + 1}\sum_{j=1}^{k}2^{-j}
	  \leq    \delta_{k + 1} \sum_{j=1}^{k + 1}2^{-j}\,,
	$$
	which is of course follows from \eqref{dafare}.
\EP

Now that we  reduced $\cL_0$ to the form $\cL_m = \omega \cdot \partial_\f + {\cal D}_m + {\cal P}_m$  
we can apply a ``standard'' KAM scheme to complete the diagonalization.
This is a super-exponentially convergent iterative scheme based on iterating the following KAM step.

\begin{lemma}[The $(m+1)$-th step]\label{passokam}
	Following the notation of Lemma \ref{lemma iterativo riducibilita} we define
	\[
	\calE_{m+1} := \Big\{\omega \in {\calE}_{m} : 
	|\omega \cdot \ell + \Omega_{m}(j) - \Omega_{ m}(j')|  \geq \frac{\gamma |j^3 - j'^3|}{\mathtt d(\ell)},  \\
	 \quad \forall (\ell , j, j') \neq (0, j, j), \quad |\ell|_\eta \leq N_{m} \Big\}
	\]
	and fix any $\zeta$ such that
	\begin{equation}\label{probab}
	 e^{- N_m\ze}   \delta_m +
  e^{C\ze^{- \mu}}  \delta_m^2\ll \de_{m+1}
	\end{equation}
Then there exists a change of variables $\Phi_{m} := {\rm exp}(\Psi_{m})$, such that
	\begin{equation}\label{psim+}
	\begin{aligned}
	\|\Psi_{m} \|_{\rho_m-\ze}^{{\calE}_{m+1}} 
	\lesssim 
	e^{C \zeta^{- \mu}} \delta_{m} \\
	\end{aligned}
	\end{equation}
	which conjugates $\cL_m$ to the operator
	\[
\cL_{m+1}= \omega \cdot \partial_\f + {\cal D}_{m+1} + {\cal P}_{m+1}\,.
	\]

The operator $\DD_{m+1}$ is of the form
\eqref{cal D (k) rid}  and satisfies \eqref{Omega j (k + 1) (k)}, with $k\rightsquigarrow m+1$,
while the operator ${\cal P}_{m+1}$ is such that
\begin{equation}\label{blabla}
\|{\cal P}_{m+1}\|_{\rho_m-\ze}^{{\calE}_{m+1}} \le  \delta_{m+1 } \,.
\end{equation}

\end{lemma}

\prova
We reason similarly to Lemma \ref{lemma iterativo riducibilita} 
i.e. we fix $\Psi_{m}$ in such a way that 
\begin{equation}\label{stoca}
\begin{aligned}
& \omega \cdot \partial_\f \Psi_m + [{\cal D}_m, \Psi_m]+ \Pi_{N_m} {\cal P}_m = {\cal Z}_m\,, \\
& {\cal Z}_m := {\rm diag}_{j \in \Z \setminus \{ 0 \}}  ({\cal P}_m)_j^j(0)\,,
\end{aligned}
\end{equation}
so that we obtains
 \begin{equation}\label{zzo}
 \|\Psi_m \|_{\rho_m - \zeta}^{{\calE}_{m+1} }
 \lesssim  e^{C \zeta^{- \mu}} \|{\cal P}_m \|_{\rho_m}^{{\calE}_{m} }
\lesssim e^{C \zeta^{- \mu}}  \delta_m\,, \\
\end{equation}
for any $\zeta<\rho_m$.

Now, for any $j \in \Z \setminus \{ 0 \}$ and any 
$\omega \in {\calE}_{m}$  one has
 $|({\cal P}_m)_j^j(0)| \lesssim \|{\cal P}_m\|_{\rho_m}^{{\calE}_{m} }
{\leq} 2\delta_m $. 
 The Hamiltonian structure guarantees that ${\cal P}_m(0)_j^j$ is purely imaginary 
 and by the Kiszbraun Theorem there exists a Lipschitz extension 
 $\omega \in {\cal O} \to \ii z_m(j)$ (with $z_m(j)$ real) 
 of this function satisfying the bound $|z_m(j)|^\calO \lesssim \delta_m$. Then, we define
	\begin{equation}\label{regalo}
	\begin{aligned}
	& {\cal D}_{m + 1} :=  {\rm diag}_{j \in \Z \setminus \{ 0 \}} \Omega_{m + 1}(j)\,, \\
	& \Omega_{m + 1}(j) := \Omega_m(j) + z_m(j) = 
	\Omega_0(j) + r_{m + 1}(j), \quad \forall j \in \Z \setminus \{ 0 \}\,, \\
	& r_{m + 1}(j) := r_m(j) + z_m(j) \\
	\end{aligned}
	\end{equation}
	and \eqref{Omega j (k + 1) (k)}, with $k\rightsquigarrow m+1$.

In order to obtain the bound \ref{blabla} we start by recalling that
	\begin{equation}\label{molta}
	 {\cal P}_{m + 1} := \Pi_{N_m}^\bot {\cal P}_m + 
	 \sum_{p \geq 2} \frac{{\rm Ad}_{\Psi_m}^{p - 1}({\cal Z}_m - \Pi_{N_m}{\cal P}_m)}{p !}  +
	 \sum_{p \geq 1} \frac{{\rm Ad}^p_{\Psi_m}({\cal P}_m)}{p!}\,, 
	\end{equation}
so that reasoning as in \eqref{stima finale cal P k + 1 (k + 1)} we obtain
	\begin{equation}\label{alledieci}
	\|{\cal P}_{m + 1}  \|_{\rho_{m}-\ze}^{\calE_{m+1}} 
	\le C' e^{- N_m\ze}   \delta_m +
  C'e^{C\ze^{- \mu}}  \delta_m^2 
	\end{equation}
and by \eqref{probab} the assertion follows.
\EP

We now iterate the step of Lemma \ref{passokam}, using at each step a smaller loss of
analyticity, namely at the $p$-th step we take $\ze_p$ with
$$
\sum_{p\ge m+1}\ze_p =\ze,
$$
so that we obtain the following standard reducibility result; for a complete proof see \cite{MP}.

	\begin{prop}\label{kam}	
	For any $j \in \Z \setminus \{ 0 \}$, the sequence $\Omega_k(j) = \Omega_0(j) + r_k(j)$, $k \geq 1$ 
	provided in Lemmata \ref{lemma iterativo riducibilita}, \ref{passokam}, 
	and defined for any $\omega \in {\cal O}$ converges to 
	$\Omega_\infty(j) = \Omega_0(j) + r_\infty(j)$ with 
	$|r_\infty(j) - r_k(j)|^{\cal O} \lesssim  \delta_k$.
        Defining the  Cantor set
	\begin{equation}\label{cantor finale cal O infty}
	\begin{aligned}
	{\calE}_\infty & := \Big\{ \omega \in {\cal O} : |\omega \cdot \ell + \Omega_\infty(j) - \Omega_\infty(j')| 
	\geq \frac{2 \gamma |j^3 - j'^3|}{\mathtt d(\ell)}, \quad \forall (\ell, j, j') \neq (0, j, j) \Big\}
	\end{aligned}
	\end{equation}
	and  
	\begin{equation}\label{def cal D infty}
	{\cal L}_\infty := \omega \cdot \partial_\f + {\cal D}_\infty , \quad {\cal D}_\infty := 
	\ii\, {\rm diag}_{j \in \Z \setminus \{  0 \}} \Omega_\infty(j)\,,
	\end{equation}
	one has 
	${\calE}_\infty \subseteq \cap_{k \geq 0} {\calE}_k$. 
	
	Defining also
	\begin{equation}\label{corona virus 10}
	\widetilde \Phi_k := \Phi_0 \circ \ldots \circ \Phi_k \quad \text{with inverse} 
	\quad \widetilde \Phi_k^{- 1} = \Phi_k^{- 1} \circ \ldots \circ \Phi_0^{- 1}\,,
	\end{equation}
	the sequence $\widetilde \Phi_k$ converges for any $\omega \in {\calE}_\infty$ to a symplectic, 
	invertible map $\Phi_\infty$ w.r.t. the norm $\| \cdot \|_{\rho_m - 2\zeta}^{{\calE}_\infty}$ and 
	$\| \Phi_\infty^{\pm 1} - {\rm Id} \|_{\rho_m - 2\zeta}^{{\calE}_\infty} \lesssim \delta_0 $. 
	Moreover for any $\omega \in {\calE}_\infty$, one has that 
	$\Phi_\infty^{- 1} {\cal L}_0 \Phi_\infty = {\cal L}_\infty$. 
\end{prop}

%%%%%%%%%%%%%%%%%%%%%%%%%%%%%%%%%%%%%%%%%%%%%%%%%
%%%%%%%%%%%%%%%%%%%%%%%%%%%%%%%%%%%%%%%%%%%%%%%%%
\subsection{Variations}
\label{meloinculo}
%%%%%%%%%%%%%%%%%%%%%%%%%%%%%%%%%%%%%%%%%%%%%%%%%
%%%%%%%%%%%%%%%%%%%%%%%%%%%%%%%%%%%%%%%%%%%%%%%%%

We now consider an operator 
\begin{equation}\label{operatore cal L0 +}
\begin{aligned}
& {\cal L}_0^+ \equiv {\cal L}_0 (\lambda_1^+, \lambda_3^+, {\cal P}_0^+) = 
\omega \cdot \partial_\f + {\cal D}_0^+ + {\cal P}_0^+\,, \\
& {\cal D}_0^+ := \lambda_3^+ \partial_x^3 + \lambda_1^+ \partial_x =
 \ii \,{\rm diag}_{j \in \Z \setminus \{ 0 \}} \Omega_0^+(j)\,, \\
& \Omega_0^+(j) := - \lambda_3^+ j^3 + \lambda_1^+ j, \quad j \in \Z \setminus \{ 0 \}\,. 
\end{aligned}
\end{equation}
such that
\begin{equation}\label{vicinanza stime riducibilita}
|\lambda_1^+ - \lambda_1|^{\calO^+}\,,\,|\lambda_3^+ - \lambda_3|^{\calO^+}\,,\,
\| {\cal P}_0^+ - {\cal P}_0 \|_{\rho_{m+1}}^{\calO^+}
 \le \delta_{m+1}
\end{equation}
where ${\cal L}$, $\lambda_1$, $\lambda_3$, ${\cal P}_0$ are given in \eqref{cal D (0)} and
${\calO^+}\subseteq\calO$. 
In other words, $\LL_0^+$ is a small variation of $\LL_0$ in \eqref{cia} with also $m\rightsquigarrow m+1$.

Of course we can apply Proposition \ref{kam} to $\LL_0^+$; our aim is to compare the ``final frequencies''
of $\LL_\io^+$ with those of $\LL_\io$.

To this aim, we first apply
 Lemma \ref{lemma iterativo riducibilita} with $\LL_0\rightsquigarrow\LL_0^+$ and
$\g\rightsquigarrow \g_+<\g$. In this way we obtain a sequence 
of sets $\calE_k^+\subseteq\calE_{k-1}^+$ and a sequence of
symplectic maps $\Phi_k^+$ defined for $\om\in\calE_{k+1}^+$ such that setting 
$\LL_0^+$ as in \eqref{operatore cal L0 +} and 
\begin{equation}\label{coniugio riducibilita}
\LL_{k}:= \Phi_{k - 1}^{- 1} {\cal L}_{k - 1} \Phi_{k - 1},
\end{equation}
one has
\begin{equation}\label{cal L (k) red}
	  {\cal L}_k^+ := \omega \cdot \partial_\f + {\cal D}_k^+ + {\cal P}_k^+,\qquad
	  k\le m+1,
\end{equation}
where
	\begin{equation}\label{cal D (k) rid +}
	\begin{aligned}
	& {\cal D}_k^+ := {\rm diag}_{j \in \Z \setminus \{ 0 \}} \Omega_k^+(j), \quad \Omega_k^+(j) 
	= \Omega_0^+ (j) + r_k^+(j)
	\end{aligned}
	\end{equation}

The sets ${\cal E}_k^+$ are defined as ${\cal E}_0^+ := {\cal O^+}$ and for $k \geq 1$
	\begin{equation}\label{def cal O (k) +}
	\begin{aligned}
		{\cal E}_k^+  := \Big\{\omega \in &{\cal E}_{k - 1}^+ \ :\  
		|\omega \cdot \ell + \Omega_{k - 1}^+(j) - \Omega_{k - 1}^+(j')| 
		 \geq \frac{\gamma_+ |j^3 - j'^3|}{\mathtt d(\ell)},  \\
		& \quad \quad \forall (\ell , j, j') \neq (0, j, j), \quad |\ell_\eta| \leq N_{k - 1} \Big\}\,. \\
	\end{aligned}
	\end{equation}

Moreover
one has $\Phi_{k - 1}^+ := {\rm exp}(\Psi_{k - 1}^+)$, with
	\begin{equation}\label{tidoila}
	 \|\Psi_{k - 1}^+ \|_{\rho_k}^{{\calE}^+_k} \lesssim 
	e^{C \Delta_{k-1}^{- \mu}} \delta_{k - 1} .
	\end{equation}

The following lemma holds. 
\begin{lemma}\label{lemma iterativo riducibilita +}
For all $k = 1, \ldots , m+1$ one has
	\begin{subequations}
	\begin{align}
	\|{\cal P}_k^+ - {\cal P}_k \|_{\rho_{k}}^{{\cal E}_k \cap {\cal E}_k^+} \leq   \delta_{m+1}, 
	\label{operator cal P (k) +}\\
   |r_{k}^+(j) - r_k(j)|^{\calO\cap \calO^+}\le \de_{m+1}
		 \label{stoffa}
		 \end{align}
	\end{subequations}
and 
	\begin{equation}\label{stima Psi k - 1 lemma +}
		\|\Psi_{k - 1}^+ - \Psi_{k - 1} \|_{\rho_k}^{{\cal E}_k \cap {\cal E}_k^+} \lesssim  
		 \delta_{m+1}, 
	\end{equation}
\end{lemma}

\prova
We procede differently for $k=1,\ldots,m$ and $k=m+1$.

For the first case we argue by induction.
Assume the  statement to hold up to some $k<m$. 
We want to prove
	\begin{equation}\label{stima Psi k - Psi k +}
	\|\Psi_k^+ - \Psi_k \|_{\rho_{k+1}}^{{\cal E}_{k + 1} \cap {\cal E}_{k + 1}^+} \le \de_{m+1}\,.  
	\end{equation}

	By Lemma \ref{lemma iterativo riducibilita}, one has for $\omega \in {\cal E}_{k + 1}^+$
	\begin{equation}\label{def matrice Psi k +}
		( \Psi_k^+)_j^{j'} (\ell) :=  \begin{cases}
\dfrac{({({\cal P}_k^+)}^{(k)})_j^{j'}(\ell)}{\ii \big(\omega \cdot \ell + \Omega_{k}^+(j) - \Omega_k^+(j') \big)}, 
\quad \forall (\ell, j, j') \neq (0, j, j), \quad |\ell|_\eta \leq N_k\,, \\
0 \quad \text{otherwise},
	\end{cases}	
	\end{equation}
and direct calculation shows that for $\omega \in {\cal E}_{k + 1} \cap {\cal E}_{k + 1}^+$, one has 
	\begin{equation}\label{corona merda 0}
\left| {(\Omega_k^+(j) - \Omega_k^+(j')) -(\Omega_k(j) - \Omega_k(j'))} \right| \le \de_{m+1}|j^3-j'^3|
	\end{equation}
	and hence  
	\begin{equation}\label{corona merda 1}
	\begin{aligned}
	| (\Psi_k^+)_j^{j'}(\ell) -  (\Psi_k)_j^{j'}(\ell)|^{ {\cal E}_{k + 1} \cap {\cal E}_{k + 1}^+} 
	&\lesssim  \delta_{m+1}  \mathtt d(\ell)^3
	|({\cal P}_k^{(k)})_j^{j'}(\ell)|^{ {\cal E}_{k + 1} \cap {\cal E}_{k + 1}^+}\\
	&+ \mathtt d(\ell)^2
	|({\cal P}_k^{(k)})_j^{j'}(\ell) - (({\cal P}_k^+)^{(k)})_j^{j'}(\ell)|^{ {\cal E}_{k + 1} \cap {\cal E}_{k + 1}^+}
	\,. 
	\end{aligned}
	\end{equation}
	
Therefore, reasoning as in \eqref{stima coefficienti matrice Psi k}--\eqref{stima Psi (k)},
one uses Lemma \ref{bound per stima di Cauchy}, the smallness condition
\eqref{dafare} and the induction estimate \eqref{operator cal P (k) +}
so that \eqref{stima Psi k - Psi k +} follows.

Now, from the definition of $r_{k+1}$ in 
\eqref{def nuova diagonale (k + 1)} it follows
\begin{equation}\label{glierre}
|r_{k+1}^+(j) - r_{k+1}(j)|^{{{\cal E}_{k+1} \cap {\cal E}_{k+1}^+}} \le \de_{m+1},
\end{equation}
and by Kiszbraun Theorem applied to $r_{k+1}^+(j) - r_{k+1}(j)$, \eqref{stoffa} holds.

The estimate of ${\cal P}_{k + 1}^+ - {\cal P}_{k + 1}$ follows by explicit computation 
the difference by using the expressions provided in \eqref{def cal P (k + 1) b},
using the induction estimates \eqref{operator cal P (k)}, \eqref{operator cal P (k) +}, 
the estimate \eqref{stima Psi k - Psi k +} and by applying Lemma \ref{variazione 0}.

For $k=m+1$ the proof can be repeated word by word, the only difference being that
$\Psi_m$ is defined in \eqref{stoca} while $\Psi_m^+$ is defined in \eqref{def matrice Psi k}
with $k=m$.
\EP

%%%%%%%%%%%%%%%%%%%%%%%%%%%%%%%%%%%%%%%%%%%%%%%%%
\subsection{Conclusion of the proof}
%%%%%%%%%%%%%%%%%%%%%%%%%%%%%%%%%%%%%%%%%%%%%%%%%

To conclude the proof of Proposition \ref{ciaraggio} we star by noting that, setting
$\calO$ appearing in \eqref{risistemo} as $\calO^{(n)}$ appearing in \eqref{alvolo}, the operator
 $L_{n+1}$
appearing in \eqref{elleka} with of course $n\rightsquigarrow n+1$ is of the form \eqref{espansione di finale}
with
\[
\begin{aligned}
&\la_3=1+A_{n+1}\,,\\
&a_1^{(k)}(\f,x) = B_{k+1}(\f,x) - B_k(\f,x) \,,\\
&a_0^{(k)}(\f,x) = C_{k+1}(\f,x) - C_k(\f,x)\,.
\end{aligned}
\]
Moreover from \eqref{taglie} we have
\[
\de_k = \s_k^{-\tau_2}e^{\tC\s_k^{-\m}}\e_k,\qquad
\rho_k = s_k - 3\s_k
\]
where $s_k$, $\s_k$ and $\e_k$ are defined in \eqref{costanti}, so that $L_{n+1}$
satisfies \eqref{risistemo} with $m=n$. Thus, fixing
\[
\ze_k=\s_k,\qquad 2\ze =\s_k,
\]
the smallness conditions \eqref{ipotesi sommabilita ei} follows by definition. Hence we can apply
Lemma \ref{riduzione ordine 1} to $L_{n+1}$ obtaining an operator of the form \eqref{coniugio cal L cal L 0}.
In particular the conjugating operator $\MM$ satisfies 
\[
\| \MM - \Id \|^{\calO}_{s_n - 3\s_n}\lesssim \s_0^{-\tau_2}e^{\tC\s_0^{-\mu}}\e_0\,.
\]

We are now in the setting of Section \ref{pelle} with
\[
\rho_k = s_k - 4\s_k,\qquad
 \de_k=\s_k^{-\tau_3}e^{2\tC\s_k^{-\frac{1}{\h}+}}\e_k
\]
for some $\tau_3>0$. A direct calculation shows that the smallness conditions 
\eqref{lostrozzo}, \eqref{dafare}, \eqref{dafare2}, \eqref{probab} are satisfied provided we choose $N_k$
appropriately, so that we can apply Proposition \ref{kam}. 

In conclusion we obtain an operator $M_{n+1} =\MM \circ \Phi_\io$ (recall that ${\cal M}$ is constructed in Lemma \ref{riduzione ordine 1}) satisfying \eqref{ukappa+}, \eqref{diago+},
where $\Omega^{(n+1)}(j) := \Omega_{\io}(j)$ and $\calE^{(n+1)} =\calE_{\io}$.
Note that in particular the functions $\Omega^{(n+1)}(j)$ turn out to be of the form \eqref{omegone+}.

Finally \eqref{tutto+} follows from Lemmata \ref{lemma variazione parte partial x} and 
\ref{lemma iterativo riducibilita +} where $\LL_+$ has the role of $L_{n+1}$ while
$\LL$ has the role of $L_n$. This means that here we are taking $m\rightsquigarrow n-1$.

\appendix

%%%%%%%%%%%%%%%%%%%%%%%%%%%%%%%%%%%%%%%%%%%%%%%%%
%%%%%%%%%%%%%%%%%%%%%%%%%%%%%%%%%%%%%%%%%%%%%%%%%
 \section{Technical Lemmata}\label{tecnici}
 \zerarcounters
%%%%%%%%%%%%%%%%%%%%%%%%%%%%%%%%%%%%%%%%%%%%%%%%%
%%%%%%%%%%%%%%%%%%%%%%%%%%%%%%%%%%%%%%%%%%%%%%%%%

 We start by recalling few results proved in \cite{MP}. Of course, as already noted in 
 \cite{MP}-Remark 2.2, all the properties holding for $\calH(\T^\infty_{\sigma + \rho}, \ell^\infty)$
 hold verbatim for $\calH(\T^\infty_{\sigma + \rho}\times\TTT_{\s+\rho}, \ell^\infty)$.
 In particular, all the estimates below hold also for the Lipschitz norms $|\cdot|_\s^\Omega$
 and $\|\cdot\|_{\s}^\Omega$. Given two Banach spaces $X$, $Y$ we denote by $\BB(X,Y)$
 the space of bounded linear operators from $X$ to $Y$.
 
\begin{prop}[\bf Torus diffeomorphism]\label{lemma diffeo inverso}
	Let $\alpha \in {\mathcal H}(\T^\infty_{\sigma + \rho}, \ell^\infty)$ be real on real. 
	Then there exists a constant $\delta \in (0, 1)$ such that if 
	$ \rho^{- 1}| \alpha |_{\sigma + \rho} \leq \delta$, then the map 
	$ \f \mapsto \f+ \alpha(\f)$ is an invertible diffeomorphism of 
	$\T^\infty_\s$ (w.r.t. the $\ell^\infty$-topology)  and its inverse is of the form 
	$\vartheta \mapsto \vartheta + \widetilde \alpha(\vartheta)$, where 
	$\widetilde \alpha \in {\mathcal H}(\T^\infty_{\sigma+\frac\rho 2}, \ell^\infty)$
	 is real on real and satisfies the estimate
	  $| \widetilde \alpha |_{\sigma+\frac{\rho}{2}} \lesssim  | \alpha |_{\sigma + \rho}$. 
\end{prop}

\begin{coro}\label{usiamoquesto}
Given $\alpha \in {\mathcal H}(\T^\infty_{\sigma + \rho}, \ell^\infty)$  as in Proposition
 \ref{lemma diffeo inverso}, the operators
\begin{align}
\label{diffeofun}
&\Phi_\alpha : {\mathcal H}(\T^\infty_{\sigma + \rho}, X) \to {\mathcal H}(\T^\infty_\sigma, X), 
\quad u(\f) \mapsto u(\f + \alpha(\f))\,,\\
& \Phi_{\tilde\alpha} : {\mathcal H}(\T^\infty_{\sigma + \frac\rho 2}, X) \to
 {\mathcal H}(\T^\infty_\sigma, X), \quad u(\vartheta) \mapsto 
 u(\vartheta + {\tilde\alpha}(\vartheta))\nonumber
\end{align}
are bounded, satisfy 
\[
\| \Phi_\alpha \|_{{\mathcal B}\Big({\mathcal H}(\T^\infty_{\sigma + \rho}, X), 
{\mathcal H}(\T^\infty_\sigma, X) \Big)},
 \| \Phi_{\tilde\alpha} \|_{{\mathcal B}\Big({\mathcal H}(\T^\infty_{\sigma + \rho}, X), 
 {\mathcal H}(\T^\infty_\sigma, X) \Big)} \leq 1\,,
 \]
and for any $\f\in \T^\infty_\s$, $u\in {\mathcal H}(\T^\infty_{\sigma + \rho}, X), 
v\in {\mathcal H}(\T^\infty_{\sigma + \frac\rho 2}, X)$ one has 
\[
\Phi_{\tilde\alpha}\circ \Phi_{\alpha} u (\f)= u(\f) \,,\quad \Phi_\alpha\circ \Phi_{\tilde\alpha} v (\f)= v(\f)\,.
\]

Moreover $\Phi$ is close to the identity in the sense that
\begin{equation}\label{cambio variabile meno identita}
\| \Phi_\alpha (u) - u \|_\sigma \lesssim \rho^{- 1} | \alpha |_\sigma | u |_{\sigma + \rho}\,. 
\end{equation}
\end{coro}

Given a function $u \in {\cal H}(\T^\infty_\sigma, X)$, we define its average on the infinite dimensional torus as 
\begin{equation}\label{media toro infinito dimensionale}
\int_{\T^\infty} u(\f)\, d \f := \lim_{N \to + \infty} \frac{1}{(2 \pi)^N} \int_{\T^N} u(\f)\, d \f_1 \ldots d \f_N\,.
\end{equation}
By Lemma 2.6 in \cite{MP}, this definition is well posed and 
$$
\int_{\T^\infty} u(\f)\, d \f  = u(0)
$$
where $u(0)$ is the zero-th Fourier coefficient of $u$.

\begin{lemma}[\bf Algebra]\label{algebra analitiche}
	One has $| u v |_\sigma \leq | u |_\sigma | v |_\sigma$ for $u,v\in \calH(\TTT^\io_\s\times\TTT_\s)$. 
\end{lemma}

\begin{lemma}[\bf Cauchy estimates]\label{stime cauchy}
	Let $u \in \calH(\TTT^\io_{\s + \rho}\times\TTT_{\s + \rho})$. Then $| \del^k u |_{\sigma} \lesssim_k \rho^{- k} 
	| u |_{\sigma + \rho}$. 
\end{lemma}

\begin{lemma}[\bf Moser composition lemma]\label{moser type lemma}
	Let $f : B_R(0) \to{ \CCC}$ be an holomorphic function defined in a neighbourhood 
	of the origin $B_R(0)$ of the complex plane ${ \CCC}$. Then the composition 
	operator $F (u) := f \circ u$ is a well defined non linear map 
	${\mathcal H}(\T^\infty_\sigma \times \T_\sigma) \to {\mathcal H}(\T^\infty_\sigma \times \T_\sigma)$ and if 
	$| u |_\sigma \leq r < R$, one has the estimate $| F(u) |_\sigma \lesssim  1 + | u |_\sigma$. 
	If $f$ has a zero of order $k$ at $0$, then for any $| u |_\sigma \leq r < R$, one gets 
	the estimate $| F(u) |_\sigma \lesssim | u |_\sigma^k$. 
\end{lemma}

For any function $u \in {\mathcal H}(\T^\infty_\sigma, X)$, given $N > 0$, 
we define the projector $\Pi_N u$ as 
\[
\Pi_N u(\f) := \sum_{|\ell|_\zia \leq N}  u(\ell) e^{\ii \ell \cdot \f}\quad 
\text{and} \quad \Pi_N^\bot u := u - \Pi_N u\,.
\]

\begin{lemma}\label{bound per stima di Cauchy}
$(i)$ Let $\rho>0$. Then
$$
\sup_{\begin{subarray}{c}
\ell \in \Z^\infty_* \\
|\ell|_\zia < \infty
\end{subarray}} \prod_{i}( 1+  \langle i \rangle^{5} | \ell_i|^{5}) e^{- \rho |\ell|_\zia} \leq  
e^{{\tau} \ln\Big(\frac{\tau}{\rho} \Big) \rho^{-\frac{1}{\zia}}}
$$
for some constant $\tau = \tau(\zia) > 0$. 

\noindent
$(ii)$ Let $\rho > 0$. Then 
\[
\sum_{\ell \in \Z^\infty_*} e^{- \rho |\ell|_\zia}  \lesssim
e^{{\tau} \ln\Big(\frac{\tau}{\rho} \Big) \rho^{-\frac{1}{\zia}}}\,,
\]
for some constant $\tau = \tau(\zia) > 0$. 

\noindent
$(iii)$ Let $\al >  0$. For $N\gg 1$ one has
\begin{equation}
\label{taglio}
\sup_{\ell\in \Z^\infty_*:\; |\ell|_{\al}<N} 
\prod_{i}(1+\jap{i}^{5} |\ell_i|^{5}) \le (1+N)^{C(\al)N^{\frac{1}{1+\al}}}
\end{equation}
for some constant $C(\al)>0$ such that $C(\al)\to\io$ as $\al\to0$.
\end{lemma}

\begin{lemma}\label{cazzodilemma}
Given $u\in \calH(\TTT^\io_\s,X)$ for $X$ some Banach space,  let $g$ be a pointwise
absolutely convergent Formal Fourier series such that
\[
|g(\ell)|_X \le  \prod_{i}( 1+  \langle i \rangle^{5} | \ell_i|^{5})^{\tau'}   |u|_X,
\]
for some $\tau'>0$. Then for any $0<\rho<\s$, then $g\in \calH(\TTT^\io_{\s-\rho},X)$ and satisfies 
\[
|g|_{\s-\rho} \le e^{{\tau} \ln\Big(\frac{\tau}{\rho} \Big) \rho^{-\frac{1}{\zia}}} |u|_\s
\]
\end{lemma}

\proof
Follows directly from Lemma \ref{bound per stima di Cauchy} and Definition \ref{anali}.
\EP

\begin{lemma}\label{palline}
Recalling \eqref{def eta (ell)} and the definition of $|\ell|_1$ in \eqref{Z inf *},
one has
\begin{equation}\label{con}
\sum_{\ell \in \Z^\infty_*} \frac{|\ell|_1^3}{\mathtt d(\ell)}<\io.
\end{equation}
\end{lemma}

\prova
First of all note that for all $\ell\in\ZZZ^*_{\io}$ one has
\[
{|\ell|_1^3} \le \prod_i (1+ \jap{i}|\ell_i|)^3,
\]
which implies
\begin{equation}\label{jean}
\frac{|\ell|_1^3}{\mathtt d(\ell)} \lesssim \frac{1}{\prod_{i}(1+\jap{i}^{2} |\ell_i|^{2})}.
\end{equation}

Then we recall that (see \cite{Bjfa})  
\[
\sum_{\ell\in\ZZZ^*_\io}\frac{1}{\prod_{i}(1+\jap{i}^{2} |\ell_i|^{2})} <\io
\]
which implies \eqref{con}. 
\EP

\begin{lemma}\label{proprieta norma sigma riducibilita vphi x}
	Let $N,\sigma, \rho > 0$, $m, m' \in \R$, 
	${\mathcal R} \in {\mathcal H}(\T^\infty_\sigma, {\mathcal B}^{\sigma, m})$, 
	${\mathcal Q} \in {\mathcal H}(\T^\infty_{\sigma + \rho}, {\mathcal B}^{\sigma + \rho, m'})$.
	\\
	$(i)$ The product operator 
	${\mathcal R} {\mathcal Q} \in {\mathcal H}(\T^\infty_\sigma, {\mathcal B}^{\sigma, m + m'})$ 
	with $\|{\mathcal R} {\mathcal Q}\|_{\sigma, m + m'} \lesssim_m  \rho^{- |m|} 
	\|{\mathcal R}\|_{\sigma, m} \|{\mathcal Q}\|_{\sigma + \rho, m'}$. 
	If ${\mathcal R}(\omega), {\mathcal Q}(\omega)$ depend on a parameter 
	$\omega \in \Omega \subseteq \Dc$, then 
	$\|{\mathcal R}{\mathcal Q}\|_{\sigma, m + m'}^\Oo \lesssim_m \rho^{- (|m| + 2)} 
	\|{\mathcal R}\|_{ \sigma, m}^\Oo \|{\mathcal Q}\|_{\sigma + \rho, m'}^\Oo$. 
	If $m = m' = 0$, one has $\|{\cal R}{\cal Q}\|_\sigma^\Omega \lesssim \|{\cal R}\|_\sigma^\Omega
	 \|{\cal Q}\|_\sigma^\Omega$. 
	
	\smallskip
	
	\noindent
	$(ii)$ The projected operator 
	$\|\Pi_N^\bot {\mathcal R}\|_{\sigma, m} \leq e^{- \rho N} \|{\mathcal R}\|_{\sigma + \rho, m}$.
	
\end{lemma}

Given two linear operators ${\cal A}, {\cal B}$, we define for any $n \geq 0$, the operator 
${\rm Ad}_{{\cal A}}^n({\cal B})$ as 
$$
{\rm Ad}_{{\cal A}}^0({\cal B}) := {\cal B}, \quad  {\rm Ad}_{{\cal A}}^{n + 1}({\cal B}) := 
[{\rm Ad}_{{\cal A}}^n({\cal B}), {\cal A}]\,,
$$
where
\[
[\BB,\calA] := \BB\calA - \calA\BB \,.
\]

By iterating the estimate $(i)$ of Lemma \ref{proprieta norma sigma riducibilita vphi x}, 
one has that for any $n \geq 1$
\begin{equation}\label{stima banale Ad n}
\|{\rm Ad}_{\cal A}^n({\cal B})\|_\sigma \leq C^{n} \|{\cal A}\|^n_\sigma \|{\cal B}\|_\sigma
\end{equation}
for some constant $C > 0$.

\begin{lemma}\label{lemma commutatore R ord 1}
	Let $0  < \ldots< \rho_n < \ldots < \rho_0$ and $0 < \ldots\ll\delta_n \ll \ldots \ll \delta_0$. 
	Assume that $\sum_{i \geq 0} \delta_i < \infty$, choose any $n\ge0$
	and let ${\cal A}$ and ${\cal B}$ be linear
	operators such that
\[
\begin{aligned}
	{\cal A} = \sum_{i=0 }^n {\cal A}_i \qquad
	{\cal B} = \sum_{i = 0}^n {\cal B}_i\qquad
	\|{\cal A}_i\|_{\rho_i, - 1}\,,\, \|{\cal B}_i\|_{\rho_i, 1} \leq \delta_i, \quad i = 0, \ldots, n\,.
\end{aligned}
\]
	Then for any  $0 < \zeta_i <  \rho_i$ the following holds.
	
	\noindent
	$(i)$ For any $k \geq 1$, one has  
	$$
	{\rm Ad}_{\cal A}^k({\cal B}) = \sum_{i = 0}^n {\cal R}_i^{(k)}
	\quad\mbox{ with }\quad \|{\cal R}_i^{(k)}\|_{\rho_i - \zeta_i} \leq C_0^k \zeta_i^{- 1}  \delta_i\,
	\quad \forall i = 0, \ldots, n
	$$

	\noindent
	$(ii)$ Let ${\cal R} := e^{- {\cal A}} {\cal B} e^{ {\cal A}} - {\cal B}$. Then 
	$$
	{\cal R} = \sum_{i = 0}^n {\cal R}_i\quad\mbox{ with }\quad
	\|{\cal R}_i\|_{\rho_i - \zeta_i} \lesssim \zeta_i^{- 1}  \delta_i\, \quad \forall i = 0, \ldots, n
	$$
\end{lemma}

\prova
	{\sc Proof of item $(i)$.} We prove the statement by induction on $k$. 
	For $k = 1$, one has that 
	$$
	\begin{aligned}
	& [{\cal B}, {\cal A}]  = \sum_{i = 0}^n {\cal R}_i^{(1)}\,, \quad  {\cal R}_i^{(1)}  := [{\cal B}_i, {\cal A}_i] + \sum_{j=0}^{i-1}\big( [{\cal B}_i,{\cal A}_j ] - [{\cal A}_i , {\cal B}_j] \big)\,.
	\end{aligned}
	$$
	Since for $j < i$ one has that $ \rho_j > \rho_i$ and so all the terms in the above sum are analytic 
	at least in the strip of width $\rho_i$. 
	By applying Lemma \ref{proprieta norma sigma riducibilita vphi x}-$(i)$
	one has for any $0 < \zeta_i < \rho_i$
	$$
	\begin{aligned}
	\|{\cal R}_i^{(1)}\|_{\rho_i - \zeta_i} \lesssim \zeta_i^{- 1} 
	\Big( \delta_i^2 + \sum_{j =0}^i \delta_i \delta_j \Big) \lesssim
	 \zeta_i^{- 1} \delta_i \sum_{j \ge0} \delta_j \lesssim \zeta_i^{- 1} \delta_i
	\end{aligned}
	$$
	for $i = 0, \ldots, n$. Now we argue by induction. 
	Assume that for some $k \geq 1$, 
	${\cal R}^{(k)} := {\rm Ad}_{\cal A}^k({\cal B}) = \sum_{i = 0}^{n} {\cal R}_i^{(k)}$, with 
	$$
	\|{\cal R}_i^{(k)}\|_{ \rho_i - \zeta_i} \leq C_0^k \zeta_i^{- 1}\delta_i, \quad i = 0, \ldots, n
	$$
	for any $0 < \zeta_i < \rho_i$. Of course this implies that for all $j<i$ one has
	$$
	\|{\cal R}_j^{(k)}\|_{ \rho_i - \zeta_i} \leq C_0^k \zeta_i^{- 1}\delta_j, \quad i = 0, \ldots, n.
	$$
	
	 By definition
	$$
	\begin{aligned}
	& {\rm Ad}_{{\cal A}}^{k + 1}({\cal B}) = [{\cal R}^{(k)}, {\cal A}] = \sum_{i = 0}^n {\cal R}_i^{(k + 1)}\,, \\
	& {\cal R}_i^{(k + 1)} :=  [{\cal R}_i^{(k)}, {\cal A}_i] + \sum_{j =0}^{i-1}
	 \big( [{\cal R}_i^{(k)},{\cal A}_j ] - [{\cal A}_i , {\cal R}_j^{(k)}] \big) \,, \quad \forall i = 0, \ldots, n. 
	\end{aligned}
	$$
	Hence by applying Lemma \ref{proprieta norma sigma riducibilita vphi x}-$(i)$ 
      and using the induction hypothesis, one obtains 
	$$
	\begin{aligned}
	& \|{\cal R}_i^{(k + 1)}\|_{\rho_i - \zeta_i} \leq 
	C\Big(  \|{\cal R}_i^{(k)}\|_{\rho_i - \zeta_i} \|{\cal A}_i\|_{\rho_i - \zeta_i} +
	 \sum_{j =0}^{i-1} \|{\cal R}_i^{(k)}\|_{\rho_i - \zeta_i} \|{\cal A}_j\|_{ \rho_i - \zeta_i} + 
	 \|{\cal R}_j^{(k)}\|_{\rho_i - \zeta_i} \|{\cal A}_i\|_{\rho_i - \zeta_i} \Big) \\
	& \leq C  \zeta_i^{- 1}C_0^k \delta_i \sum_{j =0}^{ i-1} \delta_j 
	\leq C C_0^k \zeta_i^{- 1} \delta_i \sum_{j \ge0} \de_j 
	\leq C_0^{k + 1} \zeta_i^{- 1} \delta_i.
	\end{aligned}
	$$
	
\medskip
	
	\noindent
	{\sc Proof of $(ii)$.} One has
	$$
	{\cal R} = e^{- {\cal A}} {\cal B} e^{\cal A} - {\cal B} = 
	\sum_{k \geq 1} \frac{{\rm Ad}^k_{\cal A}({\cal B})}{k!} \stackrel{(i)}{=} \sum_{i = 0}^n {\cal R}_i
	\qquad\mbox{where}\qquad
	{\cal R}_i = \sum_{k \geq1} \frac{{\cal R}_i^{(k)}}{k!}\,,
	$$
	so that 
	$$
	\|{\cal R}_i\|_{\rho_i - \zeta_i} \leq \sum_{k \geq 1} \frac{\|{\cal R}_i^{(k)}\|_{\rho_i - \zeta_i}}{k!} 
	\leq \sum_{k \geq 1} \frac{C_0^k}{k!} \zeta_i^{- 1} \delta_i \lesssim \zeta_i^{- 1}\delta_i\,. 
	$$
	Therefore the assertion follows.
\EP

\begin{lemma}\label{commutatore epsilon k sigma k}
	Let $\{\rho_n\}_{n\ge0}$ and $\{\de_n\}_{n\ge0}$ as in Lemma \ref{lemma commutatore R ord 1}.
	Choose any $n>0$ and consider
	$$
	g(\f, x) = \sum_{i = 0}^n g_i(\f, x),\qquad \mbox{with }\quad
	g_i \in {\cal H}_{\rho_i},\quad | g_i |_{\rho_i} \leq \delta_i,
	\quad i=0,\ldots,n.
	$$
 Then the following holds.
	
	\noindent
	$(i)$ Consider the commutator $ [\del_x^3, {\cal G}]$ where 
	${\cal G} := \pi_0^\bot g(\f, x) \partial_x^{- 1}  $. Then, one has  
	$$
	 [\del_x^3, {\cal G}]= 3 g_x \partial_x + {\cal R}, \quad {\cal R} := \sum_{k = 0}^n {\cal R}_i,
	\qquad \mbox{where }\quad
	\|{\cal R}_i\|_{\rho_i - \zeta_i} \lesssim \zeta_i^{- 3}\delta_i,\quad\mbox{for } \  0 < \zeta_i < \rho_i.
	$$
	
	\noindent
	$(ii)$ Let $\zeta_0, \zeta_1, \ldots, \zeta_n$ satisfying $0 < 2 \zeta_i < \rho_i$, 
	$0 < \rho_n - \zeta_n < \rho_{n - 1} - \zeta_{n - 1} < \ldots < \rho_0 - \zeta_0$ and assume that 
	$\sum_{i \ge0} \zeta_i^{- 3} \delta_i < \infty$. Then, one has  
	$$
	e^{- {\cal G}}  \partial_x^3  e^{\cal G} = \partial_x^3 + 3 g_x \partial_x + {\cal R} ,
\qquad {\cal R} = \sum_{i = 0}^n {\cal R}_i ,\quad 
	\|{\cal R}_i\|_{\rho_i - 2 \zeta_i } \lesssim  \zeta_i^{- 4} \delta_i ,  \quad i = 0, \ldots, n\,.  
	$$
	
	\noindent
	$(iii)$ Let $\zeta_0, \zeta_1, \ldots, \zeta_n$ satisfying $0 < \zeta_i < \rho_i$, 
	$0 < \rho_n - \zeta_n < \rho_{n - 1} - \zeta_{n - 1} < \ldots < \rho_0 - \zeta_0$ and assume 
	that $\sum_{i \ge0} \zeta_i^{- 1} \delta_i < \infty$. Then 
	$$
	e^{- {\cal G}}  (\omega \cdot \partial_\f)   e^{\cal G} = \omega \cdot \partial_\f + {\cal R},
\qquad{\cal R} = \sum_{i = 0}^n {\cal R}_i,\quad
	\|{\cal R}_i\|_{\rho_i - \zeta_i } \lesssim \zeta_i^{- 1} \delta_i, \quad i = 0, \ldots, n\,. 
	$$
\end{lemma}

\prova
	{\sc Proof of $(i)$.} One has  
	$$
	\begin{aligned}
	& [\partial_x^3,  \pi_0^\bot g \partial_x^{- 1}]  = \pi_0^\bot (3 g_x \del_x +  3 g_{xx} + g_{xxx}\del_x^{-1}) = 
	3 g_x \partial_x + {\cal R},   \\
	& {\cal R} := \sum_{i = 0}^n {\cal R}_i, \quad {\cal R}_i := \pi_0^\bot (3 (g_i)_{xx} + (g_i)_{xxx}\del_x^{-1}) - 
	3\pi_0 (g_i)_x \partial_x\,.
	\end{aligned}
	$$
	Therefore 
	$$
	\|{\cal R}_i\|_{\rho_i -\zeta_i} \lesssim \zeta_i^{- 3} \delta_i\,. 
	$$
	
	\medskip
	
	{\sc Proof of $(ii)$.} In view of the item $(i)$, it is enough to estimate 
	$$
	\sum_{k \geq 2} \frac{{\rm Ad}_{\cal G}^k(\partial_x^3)}{k !}.
	$$
	 Let 
	\begin{equation}\label{esp Bi Ai lemma secondo}
	\begin{aligned}
	& {\cal B} := [\del_x^3, {\cal G}] = 3 g_x \partial_x + {\cal R} = \sum_{i = 0}^n {\cal B}_i, \quad 
	{\cal B}_i := 3 (g_i)_x \partial_x + {\cal R}_i\,, \quad i = 0, \ldots, n\,, \\
	& {\cal G} = \sum_{i = 0}^n {\cal G}_i, \quad {\cal G}_i := \pi_0^\bot g_i(\f, x) \partial_x^{- 1}\, 
	\quad i = 0, \ldots, n\,. 
	\end{aligned}
	\end{equation}
	One has  
	\begin{equation}\label{def delta i cal Bi}
	\begin{aligned}
	& \|{\cal B}_i\|_{\rho_i - \zeta_i, 1} \lesssim  \zeta_i^{- 3} \delta_i, \quad  \quad i = 0, \ldots, n, \\
	& \|{\cal G}_i\|_{ \rho_i - \zeta_i , - 1} \leq \|{\cal G}_i\|_{\rho_i  , - 1} 
	\lesssim | f_i |_{\rho_i} \lesssim \delta_i \leq  \zeta_i^{- 3}\delta_i, \quad i = 0, \ldots, n 
	\end{aligned}
	\end{equation}

	For any $k \geq 2$ one has
	$$
	{\rm Ad}^k_{\cal G}(\partial_x^3) = {\rm Ad}^{k - 1}_{\cal G}([\partial_x^3, {\cal G}]) 
	= {\rm Ad}_{\cal G}^{k - 1}({\cal B})\,,
	$$
hence, we can apply Lemma \ref{lemma commutatore R ord 1}
 (replacing $\rho_i$ with $\rho_i - \zeta_i$ and $\delta_i$ with $ \zeta_i^{- 3}\delta_i$) obtaining  
	$$
	{\rm Ad}^k_{\cal G}(\partial_x^3) = \sum_{i = 0}^n {\cal R}_i^{(k)}
	$$
	where ${\cal R}_i^{(k)}$ satisfies 
	\begin{equation}\label{pappagallo 0}
	\|{\cal R}_i^{(k)}\|_{\rho_i - 2 \zeta_i } \leq C_0^k \zeta_i^{- 4} \delta_i, \quad i = 0, \ldots, n
	\end{equation}
	and hence by setting
	$$
	{\cal R} = \sum_{k \geq 2} \frac{{\rm Ad}_{\cal G}^k(\partial_x^3)}{k !} = \sum_{i = 0}^n {\cal R}_i
	$$
 item $(ii)$ follows. 
	
	\noindent
	{\sc Proof of item $(iii)$.} The proof can be done arguing as in the item $(ii)$, using that 
	$$
	e^{- {\cal G}}  (\omega \cdot \partial_\f)  e^{{\cal G}} 
	= \omega \cdot \partial_\f + \sum_{k \geq 1}
	 \frac{{\rm Ad}^{k-1}_{{\cal G}} (\omega \cdot \partial_\f {\cal G})}{k!}\,,\qquad
	 \mbox{where}\quad  (\omega \cdot \partial_\f {\cal G}):= 
	  \pi_0^\bot \omega \cdot \partial_\f g(\f, x) \partial_x^{- 1} \,.
	$$
\EP

{\begin{lemma}\label{variazione 0}
		Let ${\cal A}, {\cal A}_+, {\cal B}, {\cal B}_+$ be bounded operators w.r.t. a norm
		$\|\cdot\|_\s$, and define
\begin{equation}\label{bla variazione 0}
		M_{\cal A} := {\rm max}\{ \|{\cal A}_+\|_\sigma, \|{\cal A}\|_\sigma\} \,,\qquad
		M_{\cal B} := {\rm max}\{ \|{\cal B}_+\|_\sigma, \|{\cal B}\|_\sigma\} . 
\end{equation} 
		Then the following holds. \\
		$(i)$ For any $k \geq 0$, one has
		$$
		\|{\rm Ad}_{{\cal A}_+}^k({\cal B}_+) - {\rm Ad}_{{\cal A}}^k({\cal B}) \|_{\sigma} 
		\leq C_*^k M_{\cal A}^k M_{\cal B} \big(\|{\cal A}_+ - {\cal A}\|_\sigma +
		 \|{\cal B}_+ - {\cal B}\|_\sigma \big)
		 $$
		 for some constant $C_* > 0$. \\
		$(ii)$ 
		$$
		\|e^{- {\cal A}_+} {\cal B}_+ e^{{\cal A}_+} - e^{- {\cal A}} {\cal B} e^{{\cal A}}\|_\sigma 
		\lesssim \|{\cal A}_+ - {\cal A}\|_\sigma + \|{\cal B}_+ - {\cal B}\|_\sigma\,.
		$$
\end{lemma}

\prova
		{\sc Proof of $(i)$.} We argue by induction. Of course the result is trivial for $k=0$.
		 Assume that the estimate holds for some $k \geq 1$. 
	       Then 
		$$
		\begin{aligned}
		 {\rm Ad}_{{\cal A}_+}^{k + 1}({\cal B}_+) - &{\rm Ad}_{{\cal A}}^{k + 1}({\cal B}) = {\rm Ad}_{{\cal A}_+}\Big( {\rm Ad}_{{\cal A}_+}^k({\cal B}_+)\Big) - {\rm Ad}_{{\cal A}}\Big({\rm Ad}_{{\cal A}}^k({\cal B}) \Big) \\
		& = {\rm Ad}_{{\cal A}_+}\Big( {\rm Ad}_{{\cal A}_+}^k({\cal B}_+) - {\rm Ad}_{{\cal A}}^k({\cal B}) \Big) - {\rm Ad}_{{\cal A}_+ - {\cal A}}\Big({\rm Ad}_{{\cal A}}^k({\cal B}) \Big)\,.
		\end{aligned}
		$$
		Hence, by the induction hypothesis, using \eqref{bla variazione 0}, \eqref{stima banale Ad n} and Lemma \ref{proprieta norma sigma riducibilita vphi x}-$(i)$, one obtains that 
		$$
		\begin{aligned}
		\| {\rm Ad}_{{\cal A}_+}^{k + 1}({\cal B}_+) -& {\rm Ad}_{{\cal A}}^{k + 1}({\cal B})\|_\sigma 
		\lesssim \|{\cal A}_+\|_\sigma
		 \| {\rm Ad}_{{\cal A}_+}^k({\cal B}_+) - {\rm Ad}_{{\cal A}}^k({\cal B}) \|_\sigma + 
		 \|{\cal A}_+ - {\cal A}\|_\sigma C^k \|{\cal A}\|_\sigma^k \|{\cal B}\|_\sigma \\
		& \lesssim C_*^k M_{\cal A}^{k + 1} M_{\cal B} \big( \|{\cal A}_+ - {\cal A}\|_\sigma + 
		\|{\cal B}_+ - {\cal B}\|_\sigma \big) + C^k M_{\cal A}^k M_{\cal B} \|{\cal A}_+ - {\cal A}\|_\sigma \\
		& \leq C_*^{k + 1} M_{\cal A}^{k + 1} M_{\cal B}\big( \|{\cal A}_+ - {\cal A}\|_\sigma +
		\|{\cal B}_+ - {\cal B}\|_\sigma \big)
		\end{aligned}
		$$
		for some $C_* > 0$ large enough.
	\\
		{\sc Proof of $(ii)$.} It follows by  item $(i)$, using that 
		$$
		e^{- {\cal A}_+} {\cal B}_+ e^{{\cal A}_+} - e^{- {\cal A}} {\cal B} e^{{\cal A}} 
		= \sum_{k \geq 0} \frac{{\rm Ad}_{{\cal A}_+}^k({\cal B}_+) - {\rm Ad}_{{\cal A}}^k({\cal B}) }{k!}\,. 
		$$
\EP

\begin{lemma}\label{variazione 1}
	Let ${\cal A}, {\cal A}^+ , {\cal B}, {\cal B}^+$ be linear operators
	satisfying 
	$$
	\|{\cal A}\|_{\rho, - 1}, \|{\cal A}^+\|_{\rho, - 1},
	 \|{\cal B}\|_{\rho, 1}, \|{\cal B}^+\|_{\rho, 1} <C_0.
	 $$
	  Then the following holds.
\\
	$(i)$ For any $k \geq 1$, 
	$$
	\|{\rm Ad}_{{\cal A}_+}^k({\cal B}_+) - {\rm Ad}_{{\cal A}}^k({\cal B})\|_{\rho - \zeta} \leq C^k \zeta^{- 1} 
	 \big(\|{\cal A}_+ - {\cal A}\|_{\rho, - 1} + \|{\cal B}_+ - {\cal B}\|_{\rho, 1} \big)
	 $$
	 for some constant $C > 0$ depending on $C_0$. 
\\
	$(ii)$ Setting ${\cal R} := e^{- {\cal A}} {\cal B} e^{\cal A} - {\cal B}$, and 
	${\cal R}_+ := e^{- {\cal A}_+} {\cal B}_+ e^{{\cal A}_+} - {\cal B}_+$, one has
	$$
	\|{\cal R} - {\cal R}_+\|_{\rho - \zeta} \lesssim \zeta^{- 1}  \big( \|{\cal A} - {\cal A}_+\|_{\rho, - 1} + 
	\|{\cal B} - {\cal B}_+\|_{\rho, 1}\big).
	$$ 
\end{lemma}

\prova
	{\sc Proof of $(i)$.}
	We first estimate ${\rm Ad}_{{\cal A}_+}({\cal B}_+) - {\rm Ad}_{{\cal A}}({\cal B})$. One has 
	$$
	{\rm Ad}_{{\cal A}_+}({\cal B}_+) - {\rm Ad}_{{\cal A}}({\cal B}) = {\rm Ad}_{{\cal A}_+}({\cal B}_+ - {\cal B}) + {\rm Ad}_{{\cal A}_+ - {\cal A}}({\cal B}). 
	$$

	By Lemma \ref{proprieta norma sigma riducibilita vphi x}-$(i)$, one has
	\begin{equation}\label{stime Ad lemma variazioni 1}
	\|{\rm Ad}_{\cal A}({\cal B})\|_{\rho-\zeta},
	 \|{\rm Ad}_{{\cal A}_+}({\cal B}_+)\|_{\rho-\zeta} \lesssim  \zeta^{- 1},
	\end{equation}
	and 
	\begin{equation}\label{stima 1 variazione 1}
	\begin{aligned}
	& \|{\rm Ad}_{{\cal A}_+}({\cal B}_+) - {\rm Ad}_{{\cal A}}({\cal B})\|_{\rho-\zeta} \lesssim 
	\zeta^{- 1}\big(  \|{\cal A} - {\cal A}_+\|_{\rho, - 1} + \|{\cal B} - {\cal B}_+\|_{\rho, 1}\big)\,. 
	\end{aligned}
	\end{equation} 
	
	In order to estimate ${\rm Ad}^{k }_{{\cal A}_+}({\cal B}_+) - {\rm Ad}^{k }_{{\cal A}}({\cal B}) = 
	 {\rm Ad}^{k - 1}_{{\cal A}_+}{\rm Ad}_{{\cal A}_+}({\cal B}_+) -
	  {\rm Ad}^{k - 1}_{{\cal A}}{\rm Ad}_{\cal A}({\cal B})$ for any $k \geq 2$, we apply 
	  Lemma \ref{variazione 0}-$(i)$ where we replace ${\cal B}_+$ with
	   ${\rm Ad}_{{\cal A}_+}({\cal B}_+)$ and ${\cal B}$ with ${\rm Ad}_{{\cal A}}({\cal B})$, 
	   together with the estimates \eqref{stime Ad lemma variazioni 1}, \eqref{stima 1 variazione 1}. 
\\
	{\sc Proof of $(ii)$.} It follows by $(i)$ using that ${\cal R}_+ - {\cal R} = 
	\sum_{k \geq 1} \frac{{\rm Ad}^k_{{\cal A}_+}({\cal B}_+) - {\rm Ad}^k_{{\cal A}}({\cal B})}{k!}$.
\EP

\begin{lemma}\label{variazione 2}
	Let $g_+, g \in {\cal H}_\rho$, ${\cal G} := \pi_0^\bot g(\f, x) \partial_x^{- 1}$, 
	${\cal G}_+ := \pi_0^\bot g_+(\f, x) \partial_x^{- 1}$. Then the following holds.
\\
	$(i)$ The operators ${\cal R} := e^{-{\cal G}} \partial_x^3 e^{\cal G} -  \partial_x^3 - 3 g_x \partial_x $,
	 ${\cal R}_+ := e^{-{\cal G}_+} \partial_x^3 e^{{\cal G}_+} -  \partial_x^3 - 3 (g_+)_x \partial_x $ 
	 satisfy $\|{\cal R}_+ - {\cal R}\|_{\rho - \zeta} \lesssim \zeta^{- \tau } | g_+ - g |_\rho$ for some 
	 constant $\tau > 0$.
\\
	$(ii)$ The operators 
	${\cal R} := e^{-{\cal G}}  \omega\cdot  \partial_\f   e^{\cal G} - \omega \cdot \partial_\f$ 
	and ${\cal R}_+ :=  e^{-{\cal G}_+}  \omega \cdot \partial_\f   e^{{\cal G}_+} - \omega \cdot \partial_\f$
	 satisfy the estimate $\|{\cal R}_+ - {\cal R}\|_{\rho - \zeta} \lesssim \zeta^{- \tau} | g_+ - g |_\rho$, 
	 for some constant $\tau > 0$.  
\end{lemma}

\prova
	We only prove the item $(i)$. The item $(ii)$ can be proved by similar arguments. We compute 
	\begin{equation}\label{duca}
	\begin{aligned}
	& {\cal B} := [\partial_x^3,  \pi_0^\bot g \partial_x^{- 1}]  = 
	\pi_0^\bot (3 g_x \partial_x +  3 g_{xx} + g_{xxx}\del_x^{-1}) = 3 g_x \partial_x + {\cal R}_{\cal B}\,, \\
	& {\cal R}_{\cal B} := \pi_0^\bot (3 g_{xx} + g_{xxx}\del_x^{-1}) - \pi_0(3 g_x \partial_x), \\
	&{\cal B}_+ :=  [\partial_x^3,  \pi_0^\bot g_+ \partial_x^{- 1}]  = 
	\pi_0^\bot (3 (g_+)_x \partial_x +  3 (g_+)_{xx} + (g_+)_{xxx}\del_x^{-1}) = 
	3 (g_+)_x \partial_x + {\cal R}_{{\cal B}_+}\,, \\
	& {\cal R}_{{\cal B}_+} := \pi_0^\bot (3 (g_+)_{xx} + (g_+)_{xxx}\del_x^{-1}) - \pi_0(3 (g_+)_x \partial_x)\,.
	\end{aligned}
	\end{equation}
	
	Hence 
	\begin{equation}\label{pelli}
	\begin{aligned}
	{\cal R}_+ - {\cal R} & = {\cal R}_{{\cal B}_+} - {\cal R}_{\cal B} +
	 \sum_{k \geq 2} \frac{{\rm Ad}_{{\cal G}_+}^k(\partial_x^3) - {\rm Ad}_{{\cal G}_+}^k(\partial_x^3)}{k!} \\
	& \stackrel{\eqref{duca}}{=} {\cal R}_{{\cal B}_+} - {\cal R}_{\cal B} + 
	\sum_{k \geq 2} \frac{{\rm Ad}_{{\cal G}_+}^{k - 1}({\cal B}_+) - {\rm Ad}_{{\cal G}}^{k - 1}({\cal B})}{k!} \,.
	\end{aligned}
	\end{equation}
	
	By a direct calculation one can show the estimates 
	\begin{equation}\label{acquario}
	\begin{aligned}
	& \|{\cal B}\|_{\rho - \zeta, 1} \lesssim \zeta^{- 3} | g |_\rho, \quad 
	\|{\cal B}_+\|_{\rho - \zeta, 1} \lesssim \zeta^{- 3} | g_+ |_\rho\,, \\
	& \|{\cal G}\|_{\rho, - 1} \lesssim | g |_\rho, \quad \|{\cal G}_+\|_{\rho, - 1} \lesssim | g_+ |_\rho\,, \\
	& \|{\cal R}_{{\cal B}_+} - {\cal R}_{\cal B} \|_{\rho - \zeta} \lesssim 
	\zeta^{- 3} | g_+ - g |_\rho, \quad 
	\|{\cal G}_+ - {\cal G}\|_{\rho, - 1} \lesssim | g_+ - g |_\rho\,.
	\end{aligned}
	\end{equation}
	
	The latter estimates, together with Lemma \ref{variazione 1}-$(i)$ allow to deduce 
	\begin{equation}\label{acquario 0}
	\|{\rm Ad}_{{\cal G}_+}^{k - 1}({\cal B}_+) - {\rm Ad}_{{\cal G}}^{k - 1}({\cal B})\|_{\rho - \zeta} 
	\leq C^k \zeta^{- \tau}, \quad \forall k \geq 2,
	\end{equation}
	for some constant $\tau > 0$. Thus \eqref{pelli}-\eqref{acquario 0} imply the desired bound. 
\EP
}


\begin{thebibliography}{99} 
%%%%%%%%%%%%%%%%%%%%%%%%%%%%%%%%%%%%%%%%%%%%%%%%%
%%%%%%%%%%%%%%%%%%%%%%%%%%%%%%%%%%%%%%%%%%%%%%%%%





	
{\small 
				



\bibitem{BBM}
P.~Baldi, M.~Berti, R.~Montalto.
\newblock \textit{{KAM} for quasi-linear and fully nonlinear forced {K}d{V}.}
\newblock {\em Math. Ann.}, 359:471--536, 2014.

				
\bibitem{BBHM}
P. Baldi, M.  Berti, E. Haus, R. Montalto,
\textit{Time quasi-periodic gravity water
waves in finite depth}, Inventiones Math. 214 (2), 739--911, 2018.






\bibitem{BMP1}
L.~Biasco, J.E. Massetti, M.~Procesi.
\newblock\textit{{A}n abstract {B}irkhoff {N}ormal {F}orm {T}heorem and exponential
  type stability of the 1d {NLS}.}
  
  
  
\bibitem{BMP}
L.~Biasco, J.E. Massetti, M.~Procesi.
\newblock \textit{Almost periodic solutions for the $1$d {NLS}}, 2019.
\newblock Preprint.



\bibitem{B1}{
J. Bourgain,
\textit{Construction of quasi-periodic solutions for Hamiltonian perturbations of linear equations and applications to nonlinear PDE},
{Internat. Math. Res. Notices}, 1994
}





\bibitem{B2}{
J. Bourgain,
\textit{Quasi-periodic solutions of {H}amiltonian perturbations of 2{D} linear {S}chr\"odinger equations},
Ann. of Math. (2), {\bf 148}, (1998)
}




\bibitem{B3}{
J. Bourgain,
\textit{Green's function estimates for lattice {S}chr\"odinger operators and applications},
Princeton University Press, 2005
}





\bibitem{Bgafa}{
J. Bourgain,
\textit{Construction of approximative and almost periodic solutions of
              perturbed linear {S}chr\"odinger and wave equations},
Geom. Funct. Anal. 6, no.2 (1996), {201--230}
}





\bibitem{Bjfa}
J.~Bourgain.
\newblock {\it On invariant tori of full dimension for 1{D} periodic {NLS}.}
\newblock { J. Funct. Anal.}, 229(1):62--94, 2005.




\bibitem{CP}{
L. Chierchia, P. Perfetti
\textit{Second order Hamiltonian equations on f $\TTT^\infty$ and almost-periodic solutions},
{J. Diff. Equations},  116, 1995
}




\bibitem{Y}  {
H. Cong, J. Liu, Y. Shi, X. Yuan, 
{\it The stability of full dimensional KAM tori for nonlinear Schrödinger equation}, 
J. Diff. Equations,  {\bf 264}, no.7, 2018}


\bibitem{CM}
L.~Corsi, R.~Montalto.
\newblock {\it Quasi-periodic solutions for the forced {K}irchhoff equation on
  {$\Bbb T^d$}.}
\newblock { Nonlinearity}, 31(11):5075--5109, 2018.


\bibitem{CFP}
L.~Corsi, R~Feola, M.~Procesi.
\newblock {\it {F}inite dimensional invariant {K}{A}{M} tori for tame vector fields.}
\newblock Transactions of the AMS, 372 (2019), no.3 1913-1983 






\bibitem{CW} Craig W., Wayne C. E., {\it Newton's method and periodic solutions
of nonlinear wave equation}, Comm. Pure  Appl. Math. 46,
1409-1498, 1993.



\bibitem{FP}
R.~Feola and M.~Procesi.
\newblock \textit{{K}{A}{M} for quasi-linear autonomous {N}{L}{S}},
\newblock preprint 2017.





\bibitem{K} S.~Kuksin , {\it Hamiltonian perturbations of infinite-dimensional linear systems with imaginary spectrum},
Funktsional Anal. i Prilozhen., 21, 22-37, 95, 1987.




\bibitem{K2}
S.~Kuksin
{\it Perturbations of quasiperiodic solutions of infinite-dimensional Hamiltonian systems}
 (Math. USSR Izvestiya 32  (1989), 39-62)
 
 \bibitem{K3}
S. Kuksin
{\it The perturbation theory for the quasiperiodic solutions of infinite-dimensional
Hamiltonian systems and its applications to the Korteweg de Vries equation}
  (Math. USSR  Sbornik 64 (1989), 397-413)





\bibitem{KP}
S.~Kuksin, J.~P{\"o}schel.
\newblock \textit{ Invariant {C}antor manifolds of quasi-periodic oscillations for a
  nonlinear {S}chr{\"o}dinger equation}.
\newblock {\em Ann. of Math. (2)}, 143(1):149--179, 1996.



\bibitem{Liu}
S. Liu
\textit{The existence of almost-periodic solutions for 1-dimensional nonlinear
Schr\"odinger equation with quasi-periodic forcing}
J.  Math. Phys. {\bf 61}, 031502 (2020)




\bibitem{MP}{
R. Montalto, M. Procesi
\textit{Linear Schr\"odinger equation with an almost periodic potential},
Preprint arXiv:1910.12300, 2019
}







\bibitem{P} J.~P\"oschel,
{\it A KAM-Theorem for some nonlinear PDEs},  Ann. Sc. Norm. Pisa, % Cl. Sci.,
23, 119-148, 1996.




\bibitem{Po}
J.~P{\"o}schel.
\newblock {\it On the construction of almost periodic solutions for a nonlinear
  {S}chr{\"o}dinger equation.}
\newblock { Ergodic Theory Dynam. Systems}, 22(5):1537--1549, 2002.


\bibitem{Rab1}
P. H. Rabinowitz. 
\textit{Periodic solutions of nonlinear hyperbolic partial differential equations},
 Comm. Pure Appl. Math., 20:145–205, 1967.
 
 \bibitem{Rab2}
P. H. Rabinowitz. \textit{Periodic solutions of nonlinear hyperbolic partial differential equations. II}
 Comm. Pure Appl. Math., 22:15–39, 1968.




\bibitem{RLZ}{
J. Rui, B. Liu, J. Zhang, \textit{Almost periodic solutions for a class of linear Schrödinger
equations with almost periodic forcing}, J. Math. Phys. {\bf 57}, 092702 (2016). 
18
}

\bibitem{RL1}{
J. Rui, B. Liu,
\textit{Almost-periodic solutions of an almost-periodically forced wave equation}, 
J. Math. Anal. Appl. {\bf 451}, no.2 (2017), 629--658
}


%\bibitem{RL2}{
%J. Rui, B. Liu,
%\textit{Invariant tori of full dimension for higher-dimensional beam equations with almost-periodic
%forcing}, Boundary Value Problems 77  (2020). https://doi.org/10.1186/s13661-020-01374-9
%}


 
\bibitem{W}  {C.E. Wayne, {\it Periodic and quasi-periodic solutions
of nonlinear wave equations via KAM theory}, Comm. Math. Phys.
127, 479-528, 1990.}


\bibitem{XG}
X. Xu, J. Geng, 
\textit{Almost periodic solutions of one dimensional Schrödinger equation with the external parameters},
 J. Dyn. Differ. Equations 25, 435–450 (2013).


		
		
	}
	
\end{thebibliography}
\end{document}